\documentclass{article}

\usepackage{amsmath, tensor, amssymb} 
\usepackage{amsthm}
\usepackage{graphicx}
\usepackage[labelfont=bf]{caption}
\usepackage{environ}%
\usepackage{tikz-cd}
\usepackage{enumerate} 
\usetikzlibrary{cd}

\usepackage[hidelinks]{hyperref}
\usepackage[noabbrev,capitalize]{cleveref}
\usepackage{aliascnt}

\usepackage{xcolor}

\makeatletter
\def\thanks#1{\protected@xdef\@thanks{\@thanks
        \protect\footnotetext{#1}}}
\makeatother

\theoremstyle{plain}

\newtheorem{theorem}{Theorem}[section]
\newaliascnt{proposition}{theorem}
\newtheorem{proposition}[proposition]{Proposition}
\aliascntresetthe{proposition}

\newaliascnt{corollary}{theorem}
\newtheorem{corollary}[corollary]{Corollary}
\aliascntresetthe{corollary}

\theoremstyle{definition}

\newaliascnt{definition}{theorem}
\newtheorem{definition}[definition]{Definition}
\aliascntresetthe{definition}

\newaliascnt{lemma}{theorem}
\newtheorem{lemma}[lemma]{Lemma}
\aliascntresetthe{lemma}
 
 \newaliascnt{remark}{theorem}
\newtheorem{remark}[remark]{Remark}
\aliascntresetthe{remark}

\newaliascnt{example}{theorem}
\newtheorem{example}[example]{Example}
\aliascntresetthe{example}

\newtheorem*{theorem*}{Theorem}
\newtheorem*{proposition*}{Proposition}

\definecolor{grey}{RGB}{128,128,128}

\DeclareMathOperator{\Ima}{Im}

\title{Holonomy groupoids of singular subalgebroids\\ in the non-integrable case\thanks{2020 AMS subject classification: Primary 22A22, 53C29, ~ Secondary 22E05, 53D17.\\
\hspace*{11pt} Keywords: Holonomy, Lie subalgebroids, local Lie groupoids}
}

\author{Nicolas Tribouillard, Marco Zambon}

\AtEndDocument{\bigskip{\footnotesize
\textsc{Nicolas Tribouillard:}  KU Leuven, Department of Mathematics, Celestijnenlaan 200B box 2400, 3001 Leuven, Belgium. Email:
nicolas.tribouillard@kuleuven.be\\
\textsc{Marco Zambon:} KU Leuven, Department of Mathematics, Celestijnenlaan 200B box 2400, 3001 Leuven, Belgium. Email:
 marco.zambon@kuleuven.be
 }}

\date{}

\begin{document}

\maketitle

\begin{abstract}
	\noindent The holonomy groupoid of a singular foliation on a manifold $M$ was constructed by Androulidakis and Skandalis around 2006. The construction was later extended to singular subalgebroids of an integrable Lie algebroid, upon a choice of integrating Lie groupoid.\\
	\noindent In this work, we extend the construction of the holonomy groupoid even further, by allowing for non-integrable Lie algebroids. To do so, we make a choice of integrating \emph{local} Lie groupoid $G$. When the Mal'cev  completion $AC(G)$ is smooth, our construction recovers the holonomy groupoid associated with the Lie groupoid $AC(G)$. Furthermore, our holonomy groupoid has an associated $*$-algebra. We present examples in which the holonomy groupoid is realized as a quotient of the transformation groupoid of a Lie group action.
\end{abstract}
	
\tableofcontents

\section{Introduction}

In classical differential geometry, a foliation is a decomposition of a smooth manifold $M$ into immersed submanifolds, called leaves. It determines an equivalence relation on $M$, namely the subset of $M\times M$ consisting of pairs of points that belong to the same leaf. This subset generally fails to be smooth. However, by taking into account the holonomy along the leaves, one obtains a smooth space, namely Winkelnkemper’s holonomy groupoid of the foliation \cite{Winkelkemper}. Androulidakis and Skandalis \cite{AndroulidakisSkandalis} later extended this construction to singular foliations, motivated by applications in noncommutative geometry. Their work subsequently inspired Androulidakis and Zambon \cite{Zambon} to generalize the construction to singular subalgebroids of an integrable Lie algebroid $A$. This construction relies heavily on the choice of a Lie groupoid integrating $A$.\\

This paper extends the construction to singular subalgebroids of \emph{arbitrary} Lie algebroids. Given a Lie algebroid $A$, there might exist no Lie groupoid integrating it, however there always exists a local Lie groupoid whose Lie algebroid is $A$. The constructions of holonomy groupoids described above bear a close resemblance to the associative completion of a local Lie groupoid: both begin with a disjoint union of smooth building blocks and then quotient by an equivalence relation that eliminates redundancy. We exploit this resemblance to construct the holonomy groupoid of a singular subalgebroid of any Lie algebroid, whether integrable or not. While the construction is quite abstract, we can exhibit concrete examples for singular subalgebroids arising from Lie algebroid morphisms.

\subsection*{Background:  holonomy groupoids}
	
Recall first that groupoids are generalizations of groups for which the multiplication is defined only on selected pairs of elements. Lie groupoids are groupoids which are smooth. For instance, given a manifold $M$, the homotopy classes (fixing endpoints) of smooth paths in $M$ yield a Lie groupoid, whose restriction to a point $p\in M$ is precisely the familiar fundamental group of $M$ at $p$. Every  Lie groupoid has an infinitesimal counterpart, called Lie algebroid (in the above example, the Lie algebroid is the tangent bundle $TM$). However, in contrast with the case of Lie algebras, not every Lie algebroid arises from a Lie groupoid.\\
	
\textbf{Regular foliations}.\quad The holonomy groupoid of a regular foliation $\mathcal{F}$ on a manifold $M$ was originally constructed by Winkelnkemper \cite{Winkelkemper} using a concrete definition: $H(\mathcal{F})$ is the quotient of the space of leafwise paths by germs of diffeomorphisms of local slices at the endpoints. The holonomy groupoid has numerous applications, one of them being the definition of the $C^*$-algebra of $\mathcal{F}$ as the $C^*$-algebra of the Lie groupoid $H(\mathcal{F})$. The holonomy groupoid can be used as a smooth replacement of the (generally ill-behaved) leaf space $M/\mathcal{F}$, and the $C^*$-algebra serves as a nice replacement of the functions on that leaf space; this is relevant in various areas of noncommutative geometry \cite{Connes}.\\
	
\textbf{Singular foliations.}\quad If we view a regular foliation $\mathcal{F}$ as a Lie algebroid, then $H(\mathcal{F})$ can also be defined as the smallest Lie groupoid integrating $\mathcal{F}$. By smallest is meant that any other Lie groupoid $\mathcal{G}$ integrating $\mathcal{F}$ admits a surjective morphism $\mathcal{G}\rightarrow H(\mathcal{F})$ covering the identity. This point of view was used by Androulidakis and Skandalis \cite{AndroulidakisSkandalis} to define the holonomy groupoid $H(\mathcal{F})$ of a singular foliation. By definition, a \textbf{singular foliation} $\mathcal{F}$ is an involutive, locally finitely generated submodule of the compactly supported vector fields $\mathfrak{X}_c(M)$. A singular version of Frobenius' theorem allows one to view $\mathcal{F}$ as a partition of $M$ into leaves of varying dimension. In general, the holonomy groupoid constructed in \cite{AndroulidakisSkandalis} is not a Lie groupoid.

\begin{theorem*}[\cite{AndroulidakisSkandalis}]
	Given a singular foliation $\mathcal{F}\subset \mathfrak{X}_c(M)$ of a smooth manifold $M$, there exists a  topological groupoid $H(\mathcal{F})\rightrightarrows M$ such that 
	\begin{itemize}
		\item the orbits of $H(\mathcal{F})$ coincide with the leaves of $\mathcal{F}$,
		\item for every Lie groupoid $\mathcal{G}\rightrightarrows M$ inducing $\mathcal{F}$, there exists a surjective morphism of groupoids $\mathcal{G}\rightarrow H(\mathcal{F})$. 
	\end{itemize}
\end{theorem*}
Debord \cite{Debord} showed that $H(\mathcal{F})$ is a Lie groupoid if and only if $\mathcal{F}$ is projective, or equivalently, if it is isomorphic to $\Gamma_c(V)$ for some vector bundle $V\rightarrow M$. Although it is not smooth in every direction, $H(\mathcal{F})$ carries a natural diffeology and is always smooth along the source fibres. Additionally, Androulidakis and Skandalis used $H(\mathcal{F})$ to generalize the construction of the $C^*$-algebra of a regular foliation to the singular case.\\

\textbf{Singular subalgebroids of an integrable Lie algebroid.}\quad In a different setting, let $\mathcal{G}\rightrightarrows M$ be a Lie groupoid, $A_\mathcal{G}$ its Lie algebroid and $A'\rightarrow M$ a Lie subalgebroid of $A_\mathcal{G}$ over the same base. Moerdijk and Mr\v{c}un \cite{MoerdijkMrcun} defined the smallest integration $H_\text{min}$ of $A'$ for which there exists a morphism of groupoids $\Phi_\text{min}\colon H_\text{min}\rightarrow \mathcal{G}$ integrating the inclusion $A'\hookrightarrow A_\mathcal{G}$.\\
	
The above two situations  can be unified in a common framework. Both singular foliations and Lie subalgebroids over the same base fit into the common framework of singular subalgebroids, introduced by Androulidakis and Zambon in \cite{AZ6}. Given a Lie algebroid $A$, a \textbf{singular subalgebroid} $\mathcal{B}$ is an involutive, locally finitely generated submodule of $\Gamma_c(A)$. A major goal is to develop a Lie theory for singular subalgebroids. Androulidakis and Zambon constructed the holonomy groupoid of $\mathcal{B}$ under the additional assumption that the underlying Lie algebroid $A$ is integrable. Unlike the foliation case, the construction is not canonical, since it requires a choice of Lie groupoid $\mathcal{G}$ integrating $A$.
\begin{theorem*}[\cite{Zambon,AndroulidakisZambon}]
	Let $\mathcal{G}$ be a Lie groupoid and $\mathcal{B}\subset \Gamma_c(A_\mathcal{G})$ a singular subalgebroid of its Lie algebroid. There exists a topological groupoid $H(\mathcal{B}, \mathcal{G})\rightrightarrows M$ together with a canonical morphism $\Phi\colon H(\mathcal{B}, \mathcal{G})\rightarrow \mathcal{G}$ covering the identity such that 
	\begin{itemize}
		\item $\mathcal{B}$ and the inclusion $\mathcal{B}\hookrightarrow \Gamma_c(A_\mathcal{G})$ are completely determined by $H(\mathcal{B}, \mathcal{G})$ and $\Phi$,
		\item if $A_\mathcal{G} = TM$ and $\mathcal{G} = M\times M$ is the pair groupoid, then $H(\mathcal{B}, \mathcal{G})$ coincides with the holonomy groupoid  of a singular foliation as  in \cite{AndroulidakisSkandalis},
		\item if $\mathcal{B}:= \Gamma_c(A')$ for some Lie subalgebroid $A'\hookrightarrow A_\mathcal{G}$ over the same base $M$, then $H(\mathcal{B}, \mathcal{G})$ and $\Phi$ coincide with $H_\text{min}$ and $\Phi_\text{min}$ in \cite{MoerdijkMrcun},
		\item the association $(\mathcal{B}, \mathcal{G})\rightsquigarrow (H(\mathcal{B}, \mathcal{G}),\Phi)$ is functorial.
	\end{itemize}
\end{theorem*}
	
We outline the construction of $H(\mathcal{B}, \mathcal{G})$ in \cite{Zambon}, which in turn is based on the construction of \cite{AndroulidakisSkandalis}.
The construction of $H(\mathcal{B}, \mathcal{G})$  starts from a collection of manifolds $\{U_i\}_{i\in I}$, each of which is equipped with a smooth map $\varphi_i\colon U_i\rightarrow \mathcal{G}$ encoding $\mathcal{B}$. The pair $(U_i, \varphi_i)$ is called a bisubmersion for $(\mathcal{B},\mathcal{G})$. The product of two bisubmersions $(U_1,\varphi_1)$ and $(U_2,\varphi_2)$ uses the multiplication $m$ of $\mathcal{G}$ and is defined as the fibre product $U_1\times_M U_2$ together with the composition of maps
\[\begin{tikzcd}[ampersand replacement=\&]
	{U_1\times_MU_2} \& {\mathcal{G}\times_M\mathcal{G}} \& \mathcal{G}.
	\arrow["{(\varphi_1,\varphi_2)}", from=1-1, to=1-2]
	\arrow["m", from=1-2, to=1-3]
\end{tikzcd}\]
This product procedure equips the disjoint union $\sqcup_{i\in I}U_i$ with a formal multiplication that does not yet satisfy the axioms of a groupoid, but the quotient $(\sqcup_{i\in I}U_i)/\sim$ does for a suitable equivalence relation $\sim$. Moreover, this equivalence relation kills all redundant information in the sense that it identifies elements in different bisubmersions whenever they carry the same information about $\mathcal{B}$.
	
\subsection*{Main results: arbitrary singular subalgebroids and local Lie groupoids}
The definition of $H(\mathcal{B},\mathcal{G})$ in \cite{Zambon} for a singular subalgebroid $\mathcal{B}\subset \Gamma_c(A)$ greatly relies on a choice of Lie groupoid $\mathcal{G}$ integrating $A$. It is well-known that such an integration might not exist, hence the current paper is dedicated to the removal of this integrability assumption on $A$. To this end, we mimic the construction for the holonomy groupoid given in \cite{Zambon}, but replace the Lie groupoid $\mathcal{G}$ by a \textbf{local Lie groupoid} $G$ with Lie algebroid $A$. Such a local integration exists for any Lie algebroid \cite{CrainicFernandes}. Moreover, for a suitable notion of morphism, the categories of local Lie groupoids and Lie algebroids are equivalent via the Lie functor \cite{CabreraMarcutSalazar}. \\
 
Recall that local Lie groupoids are defined just as Lie groupoids, except that the multiplication $m\colon G_m \rightarrow G$ might only be defined on a neighbourhood $G_m\subset G\tensor[_s]{\times}{_t}G$ of the units and associativity of $m$ might hold only nearby the units. There is a natural functor from the category of local Lie groupoids to topological groupoids, called the associative completion: 
\begin{equation*}
	AC(G):= \left(\bigsqcup_{n>0}\underbrace{G\tensor[_s]{\times}{_t}\cdots \tensor[_s]{\times}{_t} G}_{n\text{ factors}}\right)/\approx.
\end{equation*}
Here, $\approx$ is the equivalence relation generated by 
\begin{equation*}
	(g_1, \ldots, g_n)\approx (g_1, \ldots, m(g_i,g_{i+1}), \ldots, g_n)
\end{equation*}
whenever $(g_i,g_{i+1})\in G_m$ for some $i = 1, \ldots, n-1$. Informally, $AC(G)$ is the free groupoid generated by the elements of $G$ divided by the relations $(g,h) = m(g,h)$ for $(g,h)\in G_m$. Note that the approach and terminology used here for the associative completion of a local Lie groupoid is taken from \cite{FernandesMichiels}, but the construction was already well-known for local Lie groups \cite{vanEstvanderLee}. \\
	
The holonomy groupoid in \cite{Zambon} and the associative completion of \cite{FernandesMichiels} are constructed using a similar pattern. They are both the disjoint union of small manifolds, divided by an equivalence relation and equipped with a formal multiplication/groupoid structure. The main idea of the current paper is to exploit this similarity in order to generalize the holonomy groupoid of \cite{Zambon} to singular subalgebroids of any Lie algebroid. Our main result is the following, which summarizes Theorem \ref{prop: groupoid structure using AC(G)}, Proposition \ref{prop: same as in Marco paper when AC(G) is smooth} and  Proposition \ref{prop: morphisms functioriality}.

\begin{theorem*}
	Given a local Lie groupoid $G$ and a singular subalgebroid $\mathcal{B}\subset \Gamma_c(A_G)$, there exists a topological groupoid $H(\mathcal{B}, G)$ and a canonical morphism $\Phi\colon H(\mathcal{B}, G)\rightarrow AC(G)$ covering the identity such that 
	\begin{itemize}
		\item if $AC(G)$ is a Lie groupoid, then $H(\mathcal{B}, G)$ is canonically isomorphic to $H(\mathcal{B}, AC(G))$ from \cite{Zambon} and $\Phi$ coincides with the canonical morphism $H(\mathcal{B}, AC(G))\rightarrow AC(G)$ of \cite{Zambon},
		\item the association $(\mathcal{B}, G)\rightsquigarrow (H(\mathcal{B}, G), \Phi)$ satisfies functorial properties.
	\end{itemize}
\end{theorem*}

Many singular subalgebroids arise from a morphism of Lie algebroids $F\colon A_1\rightarrow A_2$ covering the identity, as the image $\mathcal{B}:= F_*(\Gamma_c(A_1))$. For example, the singular foliation of a Lie algebroid arises as the image of its anchor map. In that case, $H(\mathcal{B}, G_2)$ can be described explicitly as a quotient of $AC(G_1)$, where $G_1$ is an integration of $A_1$. More precisely, the following proposition paraphrases \cref{prop: H is powers of K divided by I,prop: morphism is bisubmersion for induced subalgebroid}.

\begin{proposition*} 
    Let $\varphi\colon G_1\rightarrow G_2$ be a morphism of local Lie groupoids covering the identity and   $\mathcal{B}= \varphi_*(\Gamma_c(A_{G_1}))$. The holonomy groupoid of $(\mathcal{B}, G_2)$ is given by
    \begin{equation*}
        H(\mathcal{B}, G_2)\cong \frac{AC(G_1)}{\text{Bker}(\varphi)}.
    \end{equation*}
\end{proposition*}

Here $\text{Bker}(\varphi)$ is a normal subgroupoid of $AC(G_1)$, which consists of elements admitting a local bisection lying in the kernel of the map $AC(G_1)\rightarrow AC(G_2)$ induced by $\varphi$. Furthermore, the canonical map $\Phi\colon H(\mathcal{B}, G_2)\rightarrow AC(G_2)$ is the natural map induced by $\varphi$. We present many examples for the above proposition in \S \ref{sec:Examples}. In several of these examples, $AC(G_1)$ is an action Lie groupoid.\\

The paper is organized as follows. In \S \ref{sect: preliminaries}, we recall the local objects used to define holonomy groupoids: local Lie groupoids \cite{FernandesMichiels}, singular subalgebroids \cite{AZ6} and bisubmersions \cite{Zambon}. Compared to \cite{Zambon}, the definition of a bisubmersion here is slightly more general. In \S \ref{sect: main construction}, we define an equivalence of bisubmersions and construct the holonomy groupoid of any singular subalgebroid. Moreover, the convolution algebra associated to the holonomy groupoid in  Androulidakis' Appendix A of \cite{Zambon} is generalized to the non-integrable case. In \S \ref{sect: morphism}, we discuss the morphisms relating different holonomy groupoids and use them in \S \ref{sec:Examples} to calculate the holonomy groupoid explicitly in concrete examples.\\
	
\noindent{\bf Notation:} 
The time-$\lambda$ flow of a vector field $X\in \mathfrak{X}(M)$ on a manifold $M$ is denoted $\Phi_\lambda^X$. The structure maps of a local Lie groupoid $G\rightrightarrows M$ are the source $s$, target $t$, inversion $i$ and multiplication $m$. We do not assume that $G$ is Hausdorff, but its source fibres are. The Lie algebroid $A_G$ of a local Lie groupoid $G$ is the vector bundle $\ker ds|_M$ with anchor $dt|_{A_G}$ and bracket induced by right-invariant vector fields. Its compactly supported sections are denoted $\Gamma_c(A_G)$. The $n$-th fibre power of $G$ with itself will be denoted
\begin{equation*}
	G^{(n)}:= \underbrace{G\tensor[_s]{\times}{_t}\cdots \tensor[_s]{\times}{_t} G}_{n\text{ factors}} = \left\{(g_1, \ldots, g_n)\in G^n\mid \forall i = 1, \ldots, n-1\colon s(g_i) = t(g_{i+1})\right\}.
\end{equation*}
Note that this is  not the same as the set of $n$-composable arrows, as we do not require $(g_i, g_{i+1})$ to be in the domain of $m$. The category of (global) Lie groupoids is strictly contained in the category of local Lie groupoids, hence we reserve the notation $\mathcal{G}$ for global Lie groupoids. \\
		
\noindent{\bf Acknowledgements:}
We thank Ruben Louis for explaining to us  \cite[Proposition D.4]{RubenCamille}, which we adapted to state \cref{prop: equivalence between morphisms and bisections}.
We also thank Iakovos Androulidakis for discussions on the $C^*$-algebra of the holonomy groupoid, which allowed us to formulate  \cref{rem:completion}.
Both authors acknowledge partial support by FWO project G0B3523N. 
M.Z. also acknowledges partial support by the FWO and FNRS under EOS project  G0I2222N, by FWO project G014726N (Belgium), and by  Methusalem grant METH/21/03 - long term structural funding of the Flemish Government.\\

\noindent{\bf AI Disclosure:} We used Copilot to find \cref{example: local Lie groupoid over S^2 times R}; the connection to the   literature described in that example is due to ourselves.

\section{Preliminaries}\label{sect: preliminaries}
The construction of the holonomy groupoid relies on two types of structure from the literature: local Lie groupoids (for the groupoid structure) and bisubmersions (for the local charts). We first recall local Lie groupoids and singular subalgebroids, also to fix notation and conventions. We then generalize the definition of bisubmersions from \cite{Zambon} to allow for local Lie groupoids.

\subsection{Local Lie groupoids}\label{sec:locLiegroids}
		
\subsubsection*{Local Lie groupoids and their Lie algebroids}
We follow Fernandes and Michiels \cite{FernandesMichiels} for the theory of local Lie groupoids with only a few changes.
 
\begin{definition}\label{def: local Lie groupoid}
	A (topological) \textbf{local groupoid} $G \rightrightarrows M$ consists of two topological spaces $G$ (the space of arrows) and $M\subset G$ (the space of objects) and four continuous maps
	\begin{gather*}
		s\colon G \rightarrow M, \qquad t\colon G\rightarrow M,\\
		i\colon G\rightarrow G\colon g\mapsto g^{-1}, \qquad m\colon G_m \rightarrow G\colon (g_1,g_2)\mapsto g_1 g_2
	\end{gather*}
	called source, target, inversion and multiplication. The domain $G_m\subset G\tensor[_s]{\times}{_t}G $ of the multiplication map is an open neighbourhood of
	\begin{equation*}
		\left(G\tensor[_s]{\times}{_t}M\right)\cup\left(M\tensor[_s]{\times}{_t} G\right)\cup \left\{(g, g^{-1})\mid g\in G\right\},
	\end{equation*}
	hence any element $g\in G$ can be multiplied with the corresponding unit and inverse. All the axioms of a groupoid must hold except associativity, which is only required to hold on a neighbourhood of 
	\begin{equation*}
		\left(G\tensor[_s]{\times}{_t}M\tensor[_s]{\times}{_t}M\right)\cup \left(M\tensor[_s]{\times}{_t}G\tensor[_s]{\times}{_t}M\right)\cup \left(M\tensor[_s]{\times}{_t}M\tensor[_s]{\times}{_t}G\right),
	\end{equation*}
	and whenever both products $(g_1g_2)g_3$ and $g_1(g_2g_3)$ are defined.
\end{definition}

In short, a local groupoid is a groupoid whose multiplication is {only defined locally} and whose associativity is only required locally. In contrast with \cite{FernandesMichiels}, we do require the inversion map to be defined everywhere on $G$.  (Fernandes and Michiels \cite{FernandesMichiels} allow $i$ to be defined on a strict subset of $G$, because it is necessary to find local Lie groupoid structures with a prescribed Maurer-Cartan form. We will not be interested in Maurer-Cartan forms here and the definition of inverse of a bisubmersion, \Cref{def: inverse of a bisubmersion}, is considerably simpler if $i$ is defined on the whole of $G$.)

\begin{definition}
	A \textbf{local Lie groupoid} is a local groupoid for which $G$ is a manifold, $M$ a closed, embedded submanifold of $G$, the four maps involved are all smooth, and $s,t,m$ are submersions.
\end{definition}

\begin{definition}\label{def: morphism of local Lie groupoids}
	A \textbf{morphism of local Lie groupoids} from $G_1$ to $G_2$ is a smooth map $F\colon G_1\rightarrow G_2$ such that $F((G_1)_m)\subset (G_2)_m$ and intertwining all the structure maps of $G_1$ and $G_2$.
\end{definition}

\begin{remark}
	\Cref{def: morphism of local Lie groupoids} is more or less the same as in \cite{FernandesMichiels}, except that we assume the inversion maps to be defined everywhere. In \cite{CabreraMarcutSalazar}, the maps defined in \Cref{def: morphism of local Lie groupoids} are called local Lie groupoids maps, while the term ``morphism" is reserved for germs of such maps, see \Cref{rem: local integrations of Lie algebroids}. We will not use germs of local Lie groupoid maps and prefer to stick with the convention of \cite{FernandesMichiels} for simplicity.
\end{remark}

For brevity, the fibre product of $G$ with itself will be denoted $G^{(2)}:= G \tensor[_s]{\times}{_t}G$ and more generally
\begin{equation*}
	G^{(n)}:= \underbrace{G\tensor[_s]{\times}{_t}\cdots \tensor[_s]{\times}{_t} G}_{n \text{ factors}}= \{(g_1, \ldots, g_n)\in G^n\mid s(g_i) = t(g_{i+1})\}.
\end{equation*}
There are two methods to obtain smaller local Lie groupoids from a local Lie groupoid $G \rightrightarrows M$. By a \textbf{shrinking} of $G$ we mean an open neighbourhood $G' \subset G$ around $M$ such that $i(G') = G'$, together with the natural restrictions of the four maps to $G'$ and domain of multiplication
\begin{gather*}
	G_m' = \left(G'\right)^{(2)}\cap m^{-1}\left(G'\right).
\end{gather*}
The second method {yields} a \textbf{restriction} of $G$. The latter is a local groupoid $G'$ with the same set of objects and arrows as $G$, but a smaller domain $G_m'\subset G_m$ for $m$. In the theory of local Lie groupoids, it is common to take restrictions and shrinkings of a local Lie groupoid $G$ until it has the desired properties.

\begin{definition}
	A local Lie groupoid $G$ is \textbf{bi-regular} if for every $(g,h)\in G_m$, the derivative of the left and right multiplication maps
	\begin{gather*}
		dL_g\colon \ker d_ht\rightarrow \ker d_{gh}t,\\
		dR_h\colon \ker d_gs\rightarrow \ker d_{gh}s
	\end{gather*}
	are isomorphisms of vector spaces. It is said to have \textbf{products connected to the axes} if for every $(g,h)\in G_m$ there exists either a path $\gamma$ starting at $s(g)$ and ending at $h$ such that $(g,\gamma(\tau))\in G_m$, or starting at $t(h)$ and ending at $g$ such that $(\gamma(\tau), h)\in G_m$ for every $\tau\in [0,1]$. Notice that in the first case, $\gamma$ necessarily lies in a $t$-fibre, while in the second, it lies in a $s$-fibre. {Geometrically, this condition ensures that every composable pair can be connected to one of the coordinate axes within $G_m$.}
\end{definition}

\begin{definition}
	Given a natural number $n\geq 3$, a local groupoid $G$ is called $\mathbf{n}$\textbf{-associative} if the following holds for all $3 \leq k \leq n$ and all $(g_1, \ldots, g_k)\in G^{(k)}$. If there exists a bracketing on $(g_1, \ldots, g_k)$ such that the resulting product $g\in G$ is defined, then for any other bracketing of $(g_1, \ldots, g_k)$, the resulting product is either also equal to $g$ or undefined. 	
\end{definition}

In other words, $G$ is $n$-associative if for all words $(g_1, \ldots, g_k)$ in $G$ of length at most $n$, their product is either not defined or unambiguous, so that we can write $g_1 \cdot \ldots \cdot g_k$. For any $n\geq 3$, there is a restriction of $G$ that is $n$-associative, but this cannot always be done for all $n$ simultaneously \cite[Proposition 2.9]{FernandesMichiels}.\\
		
The Lie algebroid $A_G$ of a bi-regular local Lie groupoid is constructed exactly as for Lie groupoids, namely $A_G = \ker ds|_M$, the anchor is $dt|_{A_G}$ and the bracket is defined using right-invariant vector fields. 
	
\begin{center}
	{\it From now on, we will assume that $G$ has been shrunk and restricted so that it is $s$-connected, bi-regular, $3$-associative and has product connected to the axes.}
\end{center}

This is justified, since shrinking/restricting a local Lie groupoid does not change its Lie algebroid.

\begin{theorem}[\cite{CrainicFernandes}]\label{thm: Lie algebroids have local integrations}
	Every Lie algebroid integrates to a local Lie groupoid.
\end{theorem}

\begin{remark}\label{rem: local integrations of Lie algebroids}
	See \cite{CabreraMarcutSalazar} and \cite{FernandesZhang} for explicit integrations of a given Lie algebroid $A$. Note that the inversion map of the spray groupoid described in \cite[Section 4.2.2]{FernandesZhang} is defined on the whole of $G$. Cabrera, Marcut and Salazar \cite{CabreraMarcutSalazar}  consider the category of local Lie groupoids and \emph{germs} of maps as in \Cref{def: morphism of local Lie groupoids}. They show that  the Lie functor, from this category to the  one of Lie algebroids, is an equivalence of categories.
\end{remark}	

\subsubsection*{The associative completion and associators}
In an attempt to globalize a local Lie groupoid $G$ to a global Lie groupoid, consider all fibre powers $\sqcup_{n>0}G^{(n)}$ of $G$. If $(g_j, g_{j+1})\in G_m$, then the elements $g := (g_1, \ldots, g_n)\in G^{(n)}$ and 
\begin{equation*}
	h := (g_1, \ldots, g_{j-1}, m(g_j, g_{j+1}), g_{j+2}, \ldots, g_n)\in G^{(n-1)}
\end{equation*}
will be called \textbf{elementary equivalent}. The element $h$ is said to be a \textbf{contraction} of $g$, while $g$ is called an \textbf{expansion} of $h$.

\begin{definition}
	The \textbf{associative completion/Mal'cev completion} of $G$ is
	\begin{equation*}
		AC(G):= \bigsqcup_{n>0}G^{(n)}/\approx,
	\end{equation*}
	where $\approx$ is the equivalence relation generated by the relation `elementary equivalent'. Two elements in the same class of $\approx$ will be called \textbf{Mal'cev equivalent}. The set $AC(G)$ is equipped with a natural groupoid structure, where multiplication is given by concatenation of tuples and inversion by $i(g_1, \ldots, g_n) = (g_n^{-1}, \ldots, g_1^{-1})$.
\end{definition}
	
\begin{remark}\label{rem: functoriality of AC}
	Note that $AC$ is a functor from the category of local Lie groupoids to topological groupoids \cite[\S 4]{FernandesMichiels}; if $F\colon G_1 \rightarrow G_2$ is a morphism of local Lie groupoids, then {the associated morphism between associative completions is}
	\begin{equation*}
		AC(F)\colon AC(G_1)\rightarrow AC(G_2)\colon [(g_1, \ldots, g_n)]\mapsto [(F(g_1), \ldots, F(g_n))].
	\end{equation*}
	Moreover, $AC$ comes with a projection morphism $AC\colon \sqcup_{n>0}G^{(n)}\rightarrow AC(G)$ that we also denote by $AC$ (this should cause no confusion in the sequel). By definition, this morphism makes the following diagram commute.
	\[\begin{tikzcd}
		{\bigsqcup_{n>0}G_1^{(n)}} & {\bigsqcup_{n>0}G_2^{(n)}} \\
		{AC(G_1)} & {AC(G_2)}
		\arrow["F", from=1-1, to=1-2]
		\arrow["AC"', from=1-1, to=2-1]
		\arrow["AC", from=1-2, to=2-2]
		\arrow["{AC(F)}"', from=2-1, to=2-2]
	\end{tikzcd}\]
	The kernel of its restriction $AC|_G$ to $G$ deserves a special name as it controls the smoothness of $AC(G)$.
\end{remark}	

\begin{definition}\label{def: associators}
	An \textbf{associator} of $G$ is a (necessarily isotropic) element $g\in G$ for which there exists a series of contractions and expansions starting at $(g)\in G^{(1)}$ and ending at $(t(g))\in M\subset G^{(1)}$. The set of all associators of $G$ is denoted $\text{Assoc}(G)$.
\end{definition}

Said differently, $g\in \text{Assoc}(G)$ if the words of length one $(g)$ and $(t(g))$ are Mal'cev equivalent. {A fundamental result in} \cite[Theorem 4.10]{FernandesMichiels} states that $AC(G)$ has the structure of a smooth Lie groupoid for which the projection $AC$ is a submersion if and only if $\text{Assoc}(G)$ is uniformly discrete. In that case, the morphism $AC|_G\colon G\rightarrow AC(G)$ is a local diffeomorphism and $A_G = A_{AC(G)}$. \\
	
The integrability obstruction of Crainic and Fernandes \cite{CrainicFernandes} states that a Lie algebroid is integrable if and only if its monodromy groups are uniformly discrete. The precise relation between monodromy groups of $A_G$ and associators of a local integration $G$ can be found in \cite[Theorem 6.19]{FernandesMichiels}. {We now give two non-trivial examples of local Lie groupoids (and their associators) whose Lie algebroids do not integrate to a global Lie groupoids.} In \Cref{example: local Lie groupoid over S^2}, the associators of $G$ are uniform, but non-discrete, while in \Cref{example: local Lie groupoid over S^2 times R}, they are discrete in every fibre, but not uniformly. By \cite[Corollary 5.4]{CrainicFernandes}, the latter failure to integrability can only occur when $A_G$ is not transitive.

\begin{example}[{\cite[Example 6.2]{FernandesMichiels}}]\label{example: local Lie groupoid over S^2}
	Consider $M = S^2\times S^2$ and the local Lie groupoid
	\begin{equation*}
		G = \left\{(y_1, y_2, x_1,x_2,a)\in M\times M\times \mathbb{R}\mid x_1\neq -y_1, x_2\neq -y_2\right\}.
	\end{equation*}
	By $-y_1$ is meant the point on the sphere antipodal to $y_1$. The source and target maps of $G$ are given by projection onto the second and first $M$-component respectively. Fix a real number $\lambda$, the multiplication map is defined by 
	\begin{equation*}
		(z,y,a)(y,x,b) = (z,x,a+b+A(z_1,y_1,x_1)+\lambda A(z_2,y_2,x_2)),
	\end{equation*}
	where $A(z_i,y_i,x_i)$ is the signed area of the geodesic triangle spanned by the three points $z_i, y_i$ and $x_i$ on the sphere. In other words, $A(z_i,y_i,x_i) = \int_\Gamma\omega_{S^2}$ where $\omega_{S^2}\in \Omega^2(S^2)$ is the canonical volume form on the $2$-sphere (with total volume $1$) and $\Gamma$ is a homotopy contracting the geodesic triangle $\Delta(z_i,y_i,x_i)$ to a point. The above multiplication is only well-defined when $z_1\neq -x_1$ and $z_2\neq -x_2$, and it is $3$-associative if we further require that
	\begin{equation*}
		-\pi< A(z_1,y_1,x_1)<\pi, \qquad -\frac{\pi}{|\lambda|}< A(z_2,y_2,x_2) < \frac{\pi}{|\lambda|}.
	\end{equation*}
	The Lie algebroid of $G$ is the Almeida-Molino Lie algebroid $A_G = TM\oplus_\omega \mathbb{R}_M$, which does not always integrate to a global Lie groupoid. Recall that $A_G$ is the vector bundle $TM\oplus (M\times\mathbb{R})\rightarrow M$ whose anchor is the projection onto $TM$ and whose bracket is twisted by the two form $\omega = \text{pr}_1^*\omega_{S^2}+\lambda\text{pr}_2^*\omega_{S^2}\in \Omega^2(M)$ as follows
	\begin{equation*}
		\left[(X,f),(Y,g)\right] := \left([X,Y], X(g)-Y(f)+\omega(X,Y)\right)
	\end{equation*}
	for all $X,Y\in \mathfrak{X}(M)$ and $f,g\in C^\infty(M)$. By \cite[Proposition 6.3]{FernandesMichiels}, $AC(G)$ is smooth if and only if $\lambda \in \mathbb{Q}$. If $\lambda$ is irrational, the associators of $G$ cannot be uniformly discrete and indeed
	\begin{equation}\label{eq: associators of Almeida-Molino}
		\text{Assoc}(G) = \Delta M\times \Lambda,
	\end{equation}
	where $\Delta M$ is the diagonal of $M\times M$ and $\Lambda\subset \mathbb{R}$ is a countable subgroup containing $\mathbb{Z}+\lambda \mathbb{Z}$.
	{\hfill $\bigtriangleup$}
\end{example}
	
\begin{example}\label{example: local Lie groupoid over S^2 times R}
	Let $I\subset\mathbb{R}$ be an interval around zero with coordinate $a$ and let $\phi$ be the real-valued function $\phi(a) = \frac{-2a}{(1+a^2)^2}$. Consider the following local Lie groupoid over $M = S^2\times I$:
	\begin{gather*}
		G = \left\{(y,x,a,b)\in S^2\times S^2\times I \times \mathbb{R}\mid y\neq -x\right\},\\
		s(y,x,a,b) = (x,a), \qquad t(y,x,a,b) = (y,a),\qquad (y,x,a,b)^{-1} = (x,y,a,-b),\\
		(z,y,a,b)(y,x,a,b') = (z,x,a,b+b'+\phi(a)A(z,y,x)),
	\end{gather*}
	where multiplication is only defined if the geodesic triangle $\Delta(z,y,x)$ on $S^2$ is uniquely defined, i.e. when $z$ and $x$ are not antipodal. The associators at a unit $(x,a)$ form a countable subgroup $\Lambda_a\subset \mathbb{R}$ containing $\phi(a)\mathbb{Z}$ that only depends on $a$. {At $a = 0$, the multiplication is associative, hence $\Lambda_0 = \{0\}$.} In \Cref{example: map onto local Lie groupoid over S2 times R}, we will need the following relation between these groups for any $a,a'\in I$ (it can be proven by a direct computation):
	\begin{equation}\label{eq: relation between associators of S2 times I}
		\phi(a)\Lambda_{a'} = \phi(a')\Lambda_a.
	\end{equation}
	We now describe the Lie algebroid of $G$, which is a special case of Example 4.3 in \cite{CrainicFernandes}. Let $\mathcal{F} = \{S^2\times\{a\}\}_{a\in I}$ be the trivial foliation of $M$, $\omega_a:= \phi(a)\omega_{S^2}$ a $2$-form on $M$ and $A:= T\mathcal{F}\oplus_{\omega_a}\mathbb{R}_M$ the central extension by $\omega_a$ with bracket
	\begin{equation*}
		\left[(X,f), (Y,g)\right] = \left(\left[X,Y\right], X(g)-Y(f)-\omega_a(X,Y)\right).
	\end{equation*}
	To see that $A$ is indeed the Lie algebroid of $G$, remark that there is a natural embedding of $G$ onto a neighbourhood of the units of the monodromy groupoid $G(A)$ of $A$:
	\begin{equation*}
		G\rightarrow G(A)\colon (y,x,a,b)\mapsto [\gamma_{y,x,a}', b],
	\end{equation*}
	where $\gamma_{y,x,a}$ is the geodesic path joining $(x,a)$ to $(y,a)$. Here we used the description for $G(A)$ given in \cite[Remark 3.3]{Crainic}. Alternatively, $A$ can be realized as the cotangent bundle of the Poisson manifold $(M, \pi)$ with $\pi:= (1+a^2)\omega_{S^2}^{-1}$.  An isomorphism of Lie algebroids between the two constructions is 
	\begin{equation*}
		T\mathcal{F}\oplus_{\omega_a}\mathbb{R}_M\rightarrow T^*M \cong \text{pr}_{S^2}^*\left(T^*S^2\right)\oplus \text{pr}_I^*\left(T^*I\right)\colon (X, b)\mapsto \left(\frac{1}{1+a^2}\omega_{S^2}(X), bda\right).
	\end{equation*}
	Note that $(M, \pi)$ is a special case of \cite[\S 6.2]{CrainicFernandes_PoissonBrackets} once we identify $M$ with an annulus inside $\mathbb{R}^3\cong \mathfrak{su}(2)^*$. Indeed, $\pi$ is the linear Poisson structure on $\mathfrak{su}(2)^*$ multiplied by the function $R\mapsto R(1+(R-1)^2)$ depending only on the radius $R = 1+a$ (beware that the letter $a$ here has a different meaning than in \cite[\S 6.2]{CrainicFernandes_PoissonBrackets}). {\hfill$\bigtriangleup$}
\end{example}

\subsubsection*{Right-invariant extensions on $G^{(n)}$}
Given a local Lie groupoid $G$ and a section $\alpha \in \Gamma(A_G)$, its right- and left-invariant extensions are the vector fields on $G = G^{(1)}$ defined by
\begin{gather*}
	\overrightarrow{\alpha}_g := dR_g(\alpha_{t(g)}),\\
	\overleftarrow{\alpha}_g := dL_g(\alpha_{s(g)}).
\end{gather*}
The proof of the following lemma is omitted, but can be found in the proof of \cite[Proposition 2.13]{Zambon}.

\begin{lemma}\label{lemma: morphisms and right-invariant extensions}
	Let $F\colon G_1\rightarrow G_2$ be a morphism of local Lie groupoids covering the identity, $\alpha_1\in \Gamma(A_{G_1})$ and $\alpha_2:= F_*(\alpha_1)\in \Gamma(A_{G_2})$. The vector field $\overrightarrow{\alpha_1}\in \mathfrak{X}(G_1)$ is $F$-related to $\overrightarrow{\alpha_2}\in \mathfrak{X}(G_2)$. Similarly, $\overleftarrow{\alpha_1}$ is $F$-related to $\overleftarrow{\alpha_2}$.
\end{lemma}

For a natural number $n>0$, the right-invariant extensions of $\alpha$ to $G^{(n)}$ are
\begin{equation}\label{eq: basis of right-invariant}
	\begin{split}
		\overrightarrow{\alpha}^{(1)}&:= \left(\overrightarrow{\alpha}, 0, \ldots, 0\right),\\
		\overrightarrow{\alpha}^{(2)}&:= \left(i_*\overrightarrow{\alpha}, \overrightarrow{\alpha}, 0, \ldots, 0\right),\\
		\overrightarrow{\alpha}^{(3)}&:= \left(0, i_*\overrightarrow{\alpha}, \overrightarrow{\alpha}, 0, \ldots, 0\right),\\
		\vdots&\\
		\overrightarrow{\alpha}^{(n)}&:= \left(0, \ldots, 0, i_*\overrightarrow{\alpha}, \overrightarrow{\alpha}\right).
	\end{split}
\end{equation}
The upper index indicates the position of $\overrightarrow{\alpha}$ and the latter is always preceded by $i_*\overrightarrow{\alpha}$ to make it a well-defined element of $\mathfrak{X}(G^{(n)})$, except for $\overrightarrow{\alpha}^{(1)}$. Similarly, the left-invariant extensions of $\alpha$ to $G^{(n)}$ are
\begin{equation*}
	\begin{split}
		\overleftarrow{\alpha}^{(1)}:=& \left(\overleftarrow{\alpha}, i_*\overleftarrow{\alpha}, 0, \ldots, 0\right),\\
		\overleftarrow{\alpha}^{(2)}:=& \left(0, \overleftarrow{\alpha}, i_*\overleftarrow{\alpha}, 0, \ldots, 0\right),\\
		\vdots&\\
		\overleftarrow{\alpha}^{(n-1)}:= &\left(0\ldots, 0, i_*\overleftarrow{\alpha}, \overleftarrow{\alpha}\right),\\
		\overleftarrow{\alpha}^{(n)}:=& \left(0, \ldots, 0, \overleftarrow{\alpha}\right).
	\end{split}
\end{equation*}
The following two lemmas show that the Mal'cev equivalence $\approx$ interacts nicely with the flow of these vector fields. They can both be stated similarly for left-invariant extensions.

\begin{lemma}\label{lemma: first lemma equivalent elements and flows}
	Let $G$ be a local Lie groupoid and take $g = (g_1, g_2)\in G^{(2)}$ and $\alpha\in \Gamma(A_G)$ arbitrary. For $\lambda\in\mathbb{R}$ {sufficiently small}, the $2$-tuple $\Phi_\lambda^{\overrightarrow{\alpha}^{(2)}}(g)$ is Mal'cev equivalent to $g$.
\end{lemma}

Recall that $\Phi_\lambda^X$ denotes the time-$\lambda$ flow of the vector field $X$. Below we repeatedly use the fact that, if $f\colon M\rightarrow N$ is a smooth map between manifolds and $X\in \mathfrak{X}(M)$ is $f$-related to $Y\in \mathfrak{X}(N)$, then integral curves of $X$ are mapped under $f$ to integral curves of $Y$. In particular, this implies the equality $f\circ \Phi_\lambda^X = \Phi_\lambda^Y\circ f$. Notice that if $g\in G_m$, then the lemma is immediate, since $\overrightarrow{\alpha}^{(2)}$ is $m$-related to the zero vector field.

\begin{proof}
	Consider the $4$-tuple $(\Phi_\lambda^{\overrightarrow{\alpha}^{(2)}}(g), g^{-1})$, namely
	\begin{equation*}
		\left( \Phi_\lambda^{ i_*\overrightarrow{\alpha}}(g_1), \Phi_\lambda^{\overrightarrow{\alpha}}(g_2), g_2^{-1}, g_1^{-1}\right).
	\end{equation*}
	It suffices to show that it is Mal'cev equivalent to a unit. Recall that $g_2$ can always be multiplied with its inverse, hence also $\Phi_\lambda^{\overrightarrow{\alpha}}(g_2)$ can be multiplied with $g_2^{-1}$ as long as $\lambda$ is small. Their product is simply $\Phi_\lambda^{\overrightarrow{\alpha}}(t(g_2))$, since $\overrightarrow{\alpha}$ is right-invariant. Similarly, $t(g_2) = s(g_1)$ can be multiplied with $g_1^{-1}$ and hence, for $\lambda$ small, the $4$-tuple in the previous equation is Mal'cev equivalent to
	\begin{equation*}
		\left(\Phi_\lambda^{ i_*\overrightarrow{\alpha}}(g_1), \Phi_\lambda^{\overrightarrow{\alpha}}(g_1^{-1})\right) = \Phi_\lambda^{\overrightarrow{\alpha}^{(2)}}(g_1, g_1^{-1}). 
	\end{equation*}
	Now apply the multiplication map to obtain that this $2$-tuple is equivalent to $\Phi_\lambda^{ 0}(g_1\cdot g_1^{-1}) = t(g_1)$.
\end{proof}

\begin{lemma}\label{lemma: second lemma equivalent elements and flows}
	Take arbitrary $g,h\in G$ and $\alpha \in \Gamma(A_G)$ and assume that $(g)\approx (h)$. For $\lambda\in \mathbb{R}$ {sufficiently small}, the elements $\Phi_\lambda^{\overrightarrow{\alpha}}(g)$ and $\Phi_\lambda^{\overrightarrow{\alpha}}(h)$ are also Mal'cev equivalent. 
\end{lemma}

\begin{proof}
	For small $\lambda$, we have $(\Phi_\lambda^{\overrightarrow{\alpha}}(g), g^{-1})\in G_m$, hence 
	\begin{equation*}
		\left(\Phi_\lambda^{\overrightarrow{\alpha}}(g), g^{-1}\right)\approx R_{g^{-1}}\left(\Phi_\lambda^{\overrightarrow{\alpha}}(g)\right) = \Phi_\lambda^{\overrightarrow{\alpha}}(t(g)),
	\end{equation*}
	where the equality holds by the right-invariance of $\overrightarrow{\alpha}$. As $t(g)=t(h)$, this shows that $(\Phi_\lambda^{\overrightarrow{\alpha}}(g), g^{-1})$ is Mal'cev equivalent to $(\Phi_\lambda^{\overrightarrow{\alpha}}(h), h^{-1})$. The result follows, because $(g)\approx(h)$. 
\end{proof}	

The following corollary will only be used in the proof of \Cref{prop: equivalence between morphisms and bisections}.
\begin{corollary}\label{cor: flows and malcev}
	Let $g,h\in G^{(n)}$ be Mal'cev equivalent and $X_1, \ldots X_p\in \mathfrak{X}(G^{(n)})$ be vector fields, each of the form $\overrightarrow{\alpha}^{(j)}$ for some $\alpha \in \Gamma(A_G)$ and $j = 1, \ldots, n$. For $\lambda_1, \ldots, \lambda_p \in \mathbb{R}$ {sufficiently small}, the composition of flows
	\begin{gather*}
		\Phi:= \Phi_{\lambda_1}^{X_1}\circ \cdots \circ\Phi_{\lambda_p}^{X_p}
	\end{gather*}
	is defined both at $g$ and $h$, and $\Phi(g)\approx\Phi(h)$. 
\end{corollary}
\begin{proof}
It suffices to prove the statement when $p=1$. If $X_1=\overrightarrow{\alpha}^{(j)}$ for $j=1$, the statement holds by Lemma \ref{lemma: second lemma equivalent elements and flows}; if $j>1$, it follows from Lemma \ref{lemma: first lemma equivalent elements and flows}.
\end{proof}

\subsection{Singular subalgebroids and their bisubmersions}

Singular subalgebroids {provide} a direct generalization of (Stefan-Sussmann) singular foliations. They were first introduced by Androulidakis and Zambon \cite{AZ6} and later studied in \cite{Zambon}. Recall that a singular foliation $\mathcal{F}$ is a submodule of the compactly supported vector fields $\mathfrak{X}_c(M)$ that is 
\begin{enumerate}
	\item locally finitely generated: for every $x\in M$, there exists a neighbourhood $U\subset M$ around $x$ and vector fields $X_1, \ldots, X_n\in \mathcal{F}$ such that for every $X\in \mathcal{F}$, there exists $f_1, \ldots, f_n\in C^\infty(U)$ satisfying
	\begin{equation*}
		\forall y\in U\colon X_y = \sum_{i = 1}^{n}f_i(y)(X_i)_y,
	\end{equation*} 
	\item involutive: the Lie bracket of two vector fields in $\mathcal{F}$ is again in $\mathcal{F}$. 
\end{enumerate}
Moreover, if $\mathcal{F}$ consists of the compactly supported sections of a distribution, then $\mathcal{F}$ is called a \textbf{regular} foliation. The pullback of $\mathcal{F}\subset \mathfrak{X}_c(M)$ along a smooth map $\varphi\colon N \rightarrow M$ is the submodule
\begin{equation*}
	\varphi^{-1}(\mathcal{F}):= \left\{X\in \mathfrak{X}_c(N)\mid \exists f_1, \ldots, f_n\in C_c^\infty(N), \exists X_1, \ldots, X_n\in \mathcal{F}\colon d\varphi(X) = \sum_i f_i( X_i\circ \varphi)\right\},
\end{equation*}
where $d\varphi$ denotes the vector bundle map from  $TN$ to the pullback vector bundle $\varphi^*TM$. Note that $\varphi^{-1}(\mathcal{F})$ is \textbf{not} a singular foliation in general, unless $\varphi$ is transverse to $\mathcal{F}$ \cite{AndroulidakisSkandalis}. For an extensive discussion on singular foliations and transversality, see \cite{CamilleRubenLeonidSingularFoliations}. Replacing vector fields by sections of an arbitrary Lie algebroid naturally leads to the following notion.
	
\begin{definition}
	Given a Lie algebroid $A$, a \textbf{singular subalgebroid} $\mathcal{B}$ is an involutive, locally finitely generated submodule of $\Gamma_c(A)$.
\end{definition}

The basic examples of singular subalgebroids are  singular foliations (the case when $A$ is a tangent bundle) and Lie subalgebras (the case when $A$ is a Lie algebra).

\begin{example}
	If $f\colon A_1\rightarrow A_2$ is a morphism of  Lie algebroids, then $\mathcal{B}:= f_*(\Gamma_c(A_1)) $ is a singular subalgebroid of $A_2$. {Using the equivalence between Lie algebroids and local Lie groupoids (up to germs of local groupoid maps, see \Cref{rem: local integrations of Lie algebroids}),} $\mathcal{B}$ can equivalently be defined as $\mathcal{B}= F_*(\Gamma_c(A_{G_1}))$ for a morphism of local Lie groupoids $F\colon G_1\rightarrow G_2$. We will say that $\mathcal{B}$ arises from a morphism of local Lie groupoids. 
	{\hfill$\bigtriangleup$}
\end{example}
Although the holonomy groupoid was historically first defined for foliations, we will mainly be interested in the class of singular subalgebroids arising from a morphism of local Lie groupoids. This class is quite large, but it does not contain all singular foliations. For example, the foliation given in \cite[Exercise 2.36]{CamilleRubenLeonidSingularFoliations} has infinitely many generators and hence cannot be the image of a morphism of Lie algebroids.
	
\begin{definition}
	Let $G$ be a local Lie groupoid and $\mathcal{B}$ a singular subalgebroid of $A_G$.  For every integer $n\geq 1$, we associate to $\mathcal{B}$ two singular foliations on $G^{(n)}$ using the vector fields introduced in \eqref{eq: basis of right-invariant}. The \textbf{right-invariant foliation} $\overrightarrow{\mathcal{B}}^{(n)}$ is the $C^\infty_c(G^{(n)})$-span of the vector fields $\overrightarrow{\alpha}^{(j)}$ for all $\alpha\in \mathcal{B}$ and $j = 1, \ldots, n$. Similarly, its \textbf{left-invariant foliation} $\overleftarrow{\mathcal{B}}^{(n)}$ is generated by $\{\overleftarrow{\alpha}^{(j)}\mid \alpha \in \mathcal{B}, j = 1, \ldots, n\}$ over $C^\infty_c(G^{(n)})$.
\end{definition}

In the sequel, a singular subalgebroid $\mathcal{B}$ will always be considered with a local integration $G$ of the underlying Lie algebroid $A_G$. We wish to define the holonomy groupoid  $H(\mathcal{B}, G)$ of a pair $(\mathcal{B}, G)$ using the groupoid structure of $G$. The local ``charts" that will be used to construct $H(\mathcal{B}, G)$ are bisubmersions.

\begin{definition}\label{def: bisubmersion}
	A \textbf{bisubmersion} $(U, \varphi)$ for $(\mathcal{B}, G)$ is a manifold $U$ together with a smooth map $\varphi\colon U \rightarrow G^{(n)}$, for some $n\ge 1$, such that
	\begin{enumerate}[I.]
		\item $s_U:= s\circ \varphi$ and $t_U:= t\circ \varphi$ are submersions, where $s$ and $t$ denote the source of the last $G$ component and the target of the first $G$ component respectively,
		\item for all $\alpha \in \mathcal{B}$ and $1 \leq j \leq n$, there exist vector fields $X,Y\in \mathfrak{X}(U)$ that are $\varphi$-related to $\overrightarrow{\alpha}^{(j)}$ and $\overleftarrow{\alpha}^{(j)}$ respectively,
		\item $\varphi^{-1}(\overrightarrow{\mathcal{B}}^{(n)}) = \Gamma_c(\ker ds_U)$ and $\varphi^{-1}(\overleftarrow{\mathcal{B}}^{(n)}) = \Gamma_c(\ker dt_U)$. 
	\end{enumerate}
	When $n=1$, we call $(U, \varphi)$ a \textbf{simple} bisubmersion.
\end{definition}

Informally, item II. ensures that all of $\mathcal{B}$ is encoded by $(U, \varphi)$, since every right- and left-invariant extension can be lifted to $U$. Item III. ensures that $(U, \varphi)$ does not contain unnecessary additional information, since all directions tangent to the source/target fibres are generated by the infinitesimal data in $\mathcal{B}$. \\

The following definition provides a canonical choice of simple bisubmersions that always exist. Recall from \cite[Remark 1.11]{Zambon} that $\alpha_1, \ldots, \alpha_n\in \mathcal{B}$ are local generators of $\mathcal{B}$ on a neighbourhood of $x\in M$ if and only if their projections $[\alpha_1], \ldots, [\alpha_n]$ to the vector space $\mathcal{B}/I_x\mathcal{B}$  form a spanning set, where  $I_x$ is the ideal of functions on $M$ vanishing at $x$.

\begin{definition}\label{def:phbisub}
	Let $x\in M$ and denote by $I_x$ the set of smooth functions on $M$ that vanish at $x$. Suppose that $\alpha_1,\ldots,  \alpha_n \in \mathcal{B}$ are sections such that their projections to $\mathcal{B}/I_x\mathcal{B}$ form a spanning set. The associated \textbf{path-holonomy bisubmersion} is given by
	\begin{equation*}
		\varphi: U \rightarrow G: (\lambda_1, \ldots, \lambda_n, y) \mapsto \Phi_1^{\sum_i \lambda_i\overrightarrow{\alpha_i}}(y).
	\end{equation*}
	Here, $\Phi_1^X$ denotes the time-$1$ flow of $X$ and $U \subset \mathbb{R}^n \times M$ is an open neighbourhood of $(0, x)$ such that the following three local conditions are satisfied. 
	\begin{itemize}
		\item The map $\varphi$ is well-defined, i.e. the flow of $\sum_{i = 1}^{n}\lambda_i\overrightarrow{\alpha_i}$ is defined around $y$.
		\item The target $t_U:= t\circ \varphi$ is a submersion, which is possible since it is a submersion at $(0, x)$. Note that the source $s_U:= s\circ \varphi$ is the projection onto $M$ and thus always a submersion.
		\item For every $\alpha\in \mathcal{B}$, there exist functions $f_i\in C^{\infty}(s_U(U))$ such that $\alpha = \sum_{i = 1}^nf_i\alpha_i$ on $s_U(U)$, i.e. $\alpha_1, \ldots, \alpha_n$ are local generators for $\mathcal{B}$ on $s_U(U)$.
	\end{itemize}
\end{definition}

The proof that a path-holonomy bisubmersion is indeed a (simple) bisubmersion is exactly the same as that of \cite[Proposition 2.18]{Zambon}. The following lemma provides a rephrasing of condition III. in Definition \ref{def: bisubmersion}. Its proof is analogous to \cite[Lemma 2.8]{Zambon}.
\begin{lemma}\label{lemma: pointwise equality is enough for 3rd condition of bisubmersion}
	Given a smooth map $\varphi\colon U\rightarrow G^{(n)}$, write 
	\begin{gather*}
		\Gamma\left(\ker ds_U\right)^{\varphi\text{-proj}}:=
		 \left\{X\in \Gamma(\ker ds_U)\mid \exists \alpha\in \mathcal{B},\exists j\in \{1, \ldots, n\}\colon \text{ $X$ is $\varphi$-related to $\overrightarrow{\alpha}^{(j)}$}\right\}
	\end{gather*}
	for the set of vector fields on $U$ that $\varphi$-project to a generator of $\overrightarrow{\mathcal{B}}^{(n)}$. It {spans the fibre} $\ker (d_us_U)$ for every $u\in U$ if and only if $\varphi^{-1}(\overrightarrow{\mathcal{B}}^{(n)}) = \Gamma_c(\ker ds_U)$.
\end{lemma}
	
Recall that we denote the inversion map of $G$ by  $i\colon G \rightarrow G$. The same symbol will be used to denote the induced ``inversion'' on $G^{(n)}$ given by $i(g_1, \ldots, g_n) = (g_n^{-1}, \ldots, g_1^{-1})$.

\begin{definition}\label{def: inverse of a bisubmersion}
	If $\varphi\colon U \mapsto G^{(n)}$ is a bisubmersion, then its \textbf{inverse} is defined as 
	\begin{equation*}
		i\circ \varphi\colon U \rightarrow G^{(n)}.
	\end{equation*}
\end{definition}

Occasionally, the inverse of a bisubmersion $U$ will be written $U^{-1}$ for brevity.

\begin{proposition}
	The inverse of a bisubmersion is a bisubmersion. 
\end{proposition}

\begin{proof}
We check that $U^{-1}$ satisfies the three conditions in \Cref{def: bisubmersion}.
	\begin{enumerate}[I.]
		\item The maps $s_{U^{-1}} = s\circ i \circ \varphi = t_U$ and $t_{U^{-1}} = s_U$ are clearly submersions. \item Take $\alpha \in \mathcal{B}$ arbitrary and let $X$ be $\varphi$-related to $-\overleftarrow{\alpha}^{(j)}$ for some $j = 1, \ldots, n$. Then $X$ is $(i\circ \varphi)$-related to $\overrightarrow{\alpha}^{n+1-j}$, because $i_*\overrightarrow{\alpha} = -\overleftarrow{\alpha}$ {and inversion reverses the order of the factors of $G^{(n)}$}.
		\item The proof of II. shows that the module $\overrightarrow{\mathcal{B}}^{(n)}$ is mapped by $i$ bijectively to  $\overleftarrow{\mathcal{B}}^{(n)}$ and therefore 
		\begin{align*}
			\left(i \circ \varphi\right)^{-1}(\overrightarrow{\mathcal{B}}^{(n)}) &= \varphi^{-1}\left(i_*\overrightarrow{\mathcal{B}}^{(n)}\right) = \varphi^{-1}(\overleftarrow{\mathcal{B}}^{(n)})\\
			& = \Gamma_c(\ker dt_U) = \Gamma_c(\ker ds_{U^{-1}}),
		\end{align*} 
	\end{enumerate}
\end{proof}

For $k = 1, \ldots, n-1$, write $m_k$ for the multiplication map $m_k\colon G^{(n)}\rightarrow G^{(n-1)}$ of the $k$-th and $(k+1)$-th component. For local Lie groupoids, $m_k$ is only defined on $G^{(k-1)}\tensor[_s]{\times}{_t} G_m \tensor[_s]{\times}{_t} G^{(n-k-1)}$. 

\begin{proposition}\label{prop: multiplication of bisubmersion is bisubmersion}
	Let $\varphi\colon U \rightarrow G^{(n)}$ be a bisubmersion. If $\varphi(U)$ lies inside the domain of definition of $m_k$ for some $k = 1, \ldots, n-1$, then $m_k\circ \varphi\colon U \rightarrow G^{(n-1)}$ is also a bisubmersion. 
\end{proposition}

\begin{proof} We check that $m_k\circ \varphi$ satisfies the three conditions  Definition \ref{def: bisubmersion}.
	\begin{enumerate}[I.]
		\item Both $\varphi$ and $m_k\circ \varphi$ have the same source and target.
		\item The right-invariant extensions $\overrightarrow{\alpha}^{(k)}$ are all $m_k$-projectable with projection
		\begin{equation}\label{eq: multiplication and right-invariant extensions}
			dm_k(\overrightarrow{\alpha}^{(j)}) = 
			\begin{cases}
				\overrightarrow{\alpha}^{(j-1)}&\text{ if } k < j-1,\\
				0 & \text{ if } k = j-1,\\
				\overrightarrow{\alpha}^{(j)} &\text{ if } k \geq j.
			\end{cases}
		\end{equation}
		Take an arbitrary $\alpha\in \mathcal{B}$ and $j = 1, \ldots, n-1$. If $j \leq k$, then any vector field $X\in \mathfrak{X}(U)$ that is $\varphi$-related to $\overrightarrow{\alpha}^{(j)}$ is also $(m_k\circ \varphi)$-related to $\overrightarrow{\alpha}^{(j)}$. If $j > k$, then any $X\in \mathfrak{X}(U)$ that is $\varphi$-related to $\overrightarrow{\alpha}^{(j+1)}$ is $(m_k\circ \varphi)$-related to $\overrightarrow{\alpha}^{(j)}$. The left-invariant extensions $\overleftarrow{\alpha}^{(j)}$ are treated similarly.
		\item Recall that $\varphi$ and $m_k\circ \varphi$ have the same source. From \Cref{lemma: pointwise equality is enough for 3rd condition of bisubmersion}, we know that $\Gamma(\ker ds_U)^{\varphi\text{-proj}}$ generates the kernel of $ds_U$ at every point and it is enough to show that $\Gamma(\ker ds_U)^{(m_k\circ\varphi)\text{-proj}}$ also does. \Cref{eq: multiplication and right-invariant extensions} gives the following inclusion
		\begin{equation*}
			\Gamma(\ker ds_U)^{\varphi\text{-proj}}\subset \Gamma(\ker ds_U)^{(m_k\circ\varphi)\text{-proj}},
		\end{equation*}
		which finishes the proof. The second equality of item III. is shown similarly using the target.
	\end{enumerate}
\end{proof}

\begin{definition}\label{def: product of bisubmersions}
	The \textbf{product} of two bisubmersions $\varphi_1\colon U_1 \rightarrow G^{(n_1)}$ and $\varphi_2\colon U_2 \rightarrow G^{(n_2)}$ is
	\begin{equation*}
		(\varphi_1 ,\varphi_2)\colon U_1 \tensor[_{s_{U_1}}]{\times}{_{t_{U_2}}} U_2\rightarrow G^{(n_1+n_2)}.
	\end{equation*}
	Occasionally, this product of bisubmersions will be denoted $U_1\times_MU_2$ for brevity.
\end{definition}

\begin{proposition}\label{prop: product of bisubmersions using AC(G)}
	The product of two bisubmersions is again a bisubmersion. 
\end{proposition}

\begin{proof}
	Write $\varphi = (\varphi_1, \varphi_2)$, $U = U_1\times_M U_2$ and $n = n_1+n_2$ for simplicity.
	\begin{enumerate}[I.]
		\item The source $s_{U} = s_{U_2}\circ\text{pr}_{U_2}$ is a composition of two submersions and similarly for the target. 
		\item Take an arbitrary $\alpha \in \mathcal{B}$ and let $X$ be $\varphi_1$-related to $\overrightarrow{\alpha}^{(1)}\in \mathfrak{X}(G^{(n_1)})$, then $(X, 0)$ is $\varphi$-related to $\overleftarrow{\alpha}^{(1)}\in \mathfrak{X}(G^{(n)})$. Similarly, for all indices $1\leq j\leq n$ except for $j = n_1+1$, one can just use the fact that $U_1$ or $U_2$ is a bisubmersion, depending on whether $j\leq n_1$ or $j>n_1+1$, to find a lift of $\overrightarrow{\alpha}^{(j)}$. The only non-trivial case is $j = n_1+1$. Let $X\in \mathfrak{X}(U_1)$ be $\varphi_1$-related to $-\overleftarrow{\alpha}^{(n_1)} \in\mathfrak{X}(G^{(n_1)}) $ and $Y \in \mathfrak{X}(U_2)$ be $\varphi_2$-related to $\overrightarrow{\alpha}^{(1)}\in \mathfrak{X}(G^{(n_2)})$ and consider the vector field $(X, Y)\in \mathfrak{X}(U)$. It is $\varphi$-related to 
		\begin{equation*}
			\overrightarrow{\alpha}^{(n_1+1)} = (0, \ldots, 0, \underbrace{-\overleftarrow{\alpha}}_{n_1}, \underbrace{\overrightarrow{\alpha}}_{n_1+1}, 0, \ldots, 0)\in \mathfrak{X}(G^{(n)}),
		\end{equation*}
		since $i_*\overrightarrow{\alpha} = -\overleftarrow{\alpha}$. The search for vector fields $\varphi$-related to $\overleftarrow{\alpha}^{(j)}$ is similar. 		
		\item The equality
		\begin{equation*}
			\varphi^{-1}\left(\overrightarrow{\mathcal{B}}^{(n)}\right) = \Gamma_c(\ker ds_U).
		\end{equation*}
		will be shown using \Cref{lemma: pointwise equality is enough for 3rd condition of bisubmersion}. The other equality in item III. of \Cref{def: bisubmersion} is shown similarly. Given an arbitrary  vector $(X, Y)\in \ker (d_{(u_1,u_2)}s_U)$ for some $(u_1, u_2)\in U$, we will show that it can be written as a sum of vectors that can be extended to vector fields in $\Gamma(\ker ds_U)^{\varphi\text{-proj}}$.\\
		
		Since $Y\in \ker ((d_{u_2}s_{U_2}))$ and $U_2$ is a bisubmersion, $Y$ can be extended to a vector field $Y'\in \mathfrak{X}(U_2)$ of the form
		\begin{equation*}
			Y' = \sum_{k}\lambda_kY'_k,
		\end{equation*}
		where $\lambda_k \in \mathbb{R}$ and $Y'_k\in\Gamma(\ker ds_{U_2})^{\varphi_2\text{-proj}}$. After reorganizing, summing and rescaling the $Y'_k$s, we can assume that the index $k$ runs from $1$ to $n$, that $\lambda_k = 1$ and that each $Y'_k$ is $\varphi_2$-related to $\overrightarrow{\alpha_k}^{(k)}$ for some $\alpha_k\in \mathcal{B}$. Since $U_1$ is a bisubmersion, there exists a vector field $X'\in \mathfrak{X}(U_1)$ that is $\varphi_1$-related to $-\overleftarrow{\alpha_1}^{(n_1)}$. The pair $(X', Y'_1)$ forms a well-defined vector field on $U$, because
		\begin{equation*}
			ds_{U_1}(X') = ds(-\overleftarrow{\alpha}_1) = dt(\overrightarrow{\alpha_1}) = dt_{U_2}(Y'_1).
		\end{equation*}
		Moreover, $(X', Y'_1)\in \Gamma(\ker ds_U)^{\varphi\text{-proj}}$, since it is $\varphi$-related to $\overrightarrow{\alpha_1}^{(n_1+1)}$. \\
		
		The map $ds_{U_1}$ at $u_1$ sends both $X$ and $X'_{u_1}$ to the same vector, namely $d_{u_2}t_{U_2}(Y)$. \Cref{lemma: pointwise equality is enough for 3rd condition of bisubmersion} for $U_1$ allows to write their difference as a linear combination
		\begin{equation*}
			X-X'_{u_1} = \sum_j \mu_j(X'_j)_{u_1}
		\end{equation*}
		of vector fields $X'_j\in \Gamma(\ker ds_{U_1})^{\varphi_1\text{-proj}}$ at $u_1$. The original vector $(X,Y)$ can be decomposed as
		\begin{equation*}
			(X, Y) = (X-X'_{u_1}, 0)+ (X'_{u_1}, Y'_{u_2}) = \sum_j\mu_j((X_j')_{u_1}, 0) + (X'_{u_1}, (Y'_1)_{u_2}) +\sum_{k = 2}^n (0, (Y'_k)_{u_2}).
		\end{equation*} 
		Each term in the right-hand side of this equation is in $\Gamma(\ker ds_U)^{\varphi\text{-proj}}$.
	\end{enumerate}
\end{proof}

\section{Construction of the holonomy groupoid}\label{sect: main construction}
This section defines the holonomy groupoid associated with an arbitrary singular subalgebroid, generalizing the construction of \cite{Zambon}. First, we define an equivalence relation on the set of bisubmersions for a pair $(\mathcal{B}, G)$ using $AC(G)$. Second, we define the holonomy groupoid $H(\mathcal{B}, G)$ as the quotient of a family of bisubmersions by this relation. Although the definition is not canonical, it does satisfy the following two important properties.

\begin{itemize}
	\item It recovers the holonomy groupoid of \cite{Zambon} when $AC(G)$ is a Lie groupoid (\Cref{prop: same as in Marco paper when AC(G) is smooth}). This is for example the case when $G$ is already a Lie groupoid, because then $AC(G) = G$.
	\item It is functorial in the sense that if $F\colon G_1\rightarrow G_2$ is a morphism of local Lie groupoids covering the identity and $\mathcal{B}$ is a singular subalgebroid of $A_{G_1}$, then there exists a canonical surjective morphism of groupoids $H(\mathcal{B}, G_1)\rightarrow H(F_*(\mathcal{B}), G_2)$, see \Cref{prop: morphisms functioriality}.
\end{itemize}

Third, we show that $H(\mathcal{B}, G)$ has an associated $*$-algebra that naturally generalizes the one defined in \cite[Appendix A]{Zambon}. This $*$-algebra will not be used in the rest of the paper.

\subsection{Equivalence of bisubmersions}
We recall that, given a singular subalgebroid $\mathcal{B}$ of an integrable Lie algebroid and a Lie groupoid $\mathcal{G}$ integrating the latter, one can canonically associate a topological groupoid $H(\mathcal{B}, \mathcal{G})$. The construction was carried out in \cite{Zambon}, generalizing the one of \cite{AndroulidakisSkandalis}. Compared to \cite{Zambon}, the main difference here when working with a local Lie groupoid $G$ is that two bisubmersions $\varphi_i\colon U_i\rightarrow G^{(n_i)}$ with $i = 1,2$ have different codomains. However, the associative completion $AC(G)$ is able to compare elements on different powers of $G$.

\begin{definition}\label{def: morphisms using AC(G)}
	Let $\varphi_i\colon U_i\rightarrow G^{(n_i)}$ be two bisubmersions for the same pair $(\mathcal{B} , G)$. A \textbf{morphism of bisubmersions} from $U_1$ to $U_2$ is a smooth map $f\colon U_1\rightarrow U_2$ making the following diagram commute.
	\[\begin{tikzcd}
		{U_1} && U_2 \\
		{G^{(n_1)}} && {G^{(n_2)}} \\
		& {AC(G)}
		\arrow["f", from=1-1, to=1-3]
		\arrow["{\varphi_1}"', from=1-1, to=2-1]
		\arrow["{\varphi_2}", from=1-3, to=2-3]
		\arrow["AC", from=2-1, to=3-2]
		\arrow["AC"', from=2-3, to=3-2]
	\end{tikzcd}\]
	An element $u_1\in U_1$ is \textbf{equivalent} to $u_2\in U_2$ if there exists a morphism of bisubmersions $f\colon U_1' \rightarrow U_2$ locally defined on an open neighbourhood $U_1'\subset U_1$ of $u_1$ such that $f(u_1) = u_2$.
\end{definition}

Clearly, the relation ``equivalent" is reflexive and transitive, but surprisingly, it is also symmetric due to \cref{prop: equivalence between morphisms and bisections}.  In other words, it is an equivalence relation on the set of bisubmersions for $(\mathcal{B} , G)$.\\
		
Morphisms of bisubmersions will allow us to identify elements on different bisubmersions that, through the map $\varphi$, carry the same information. For example, suppose that two bisubmersions $U_1$ and $U_2$ have a non-empty intersection, then the inclusion maps $U_1\cap U_2\rightarrow U_i$ can glue $U_1$ and $U_2$ along this intersection. As another example, suppose that $\varphi\colon U\rightarrow G^{(2)}$ is a bisubmersion and that $\varphi(u)$ is in the domain of the multiplication map $m\colon G_m\rightarrow G$ for some $u\in U$. Since $G_m$ is open, there exists a neighbourhood $U'$ of $u$ on which $m\circ \varphi$ is defined and the identity map from $U'$ to itself  will be a morphism of bisubmersions between $(U', \varphi)$ and $(U',m\circ \varphi)$.

\begin{definition}
	A \textbf{local bisection} for a bisubmersion $(U, \varphi)$ is a locally closed submanifold $\mathbf{b}$ of $U$ such that $s_U|_\mathbf{b} $ and $t_U|_\mathbf{b}$ are diffeomorphisms onto open subsets of $M$.
\end{definition}

Equivalently, a local bisection may be viewed as a map $\mathbf{b}: N \rightarrow U$ defined on an open $N \subset M$ that is right inverse to $s_U$ and such that $t_U \circ \mathbf{b}$ is a diffeomorphism onto its image. The following proposition shows that morphisms of bisubmersions are locally invertible. Its statement and proof are inspired by Proposition D.4 of \cite{RubenCamille}.

\begin{proposition}\label{prop: equivalence between morphisms and bisections}
	Let $\varphi_i\colon U_i\rightarrow G^{(n_i)}$ with $i = 1,2$ be two bisubmersions for $(\mathcal{B}, G)$ and $u_i\in U_i$. The following are equivalent.
	\begin{enumerate}
		\item There exists a local morphism of bisubmersions mapping $u_1$ to $u_2$.
		\item We have $\varphi_1(u_1) \approx\varphi_2(u_2)$ and there exist local bisections $u_i\in \mathbf{b}_i\subset U_i$ s.t.
		\begin{equation*}
			AC\circ\varphi_1(\mathbf{b}_1) = AC\circ\varphi_2(\mathbf{b}_2). 
		\end{equation*}
		\item There exists a manifold $W$, an element $w\in W$ and two submersions $p_i\colon W\rightarrow U_i$ such that $p_i(w) = u_i$ and the following diagram commutes.
		\[\begin{tikzcd}
			& W \\
			{U_1} && {U_2} \\
			{G^{(n_1)}} & {AC(G)} & {G^{(n_2)}}
			\arrow["{p_1}"', from=1-2, to=2-1]
			\arrow["{p_2}", from=1-2, to=2-3]
			\arrow["{\varphi_1}", from=2-1, to=3-1]
			\arrow["{\varphi_2}"', from=2-3, to=3-3]
			\arrow[from=3-1, to=3-2]
			\arrow[from=3-3, to=3-2]
		\end{tikzcd}\]
	\end{enumerate}
\end{proposition}

\begin{proof}
	The implication from \textit{1.} to \textit{2.} is obvious and \textit{3.} to \textit{1.} follows by taking a local section of $p_1$. We now prove that \textit{2.} implies \textit{3.}. \\
			
	Proposition 4.6 of \cite{FernandesMichiels} says that Mal'cev-equivalent elements admit a common refinement in a sufficiently high fibre power $G^{(p)}$. Since $\varphi_1(u_1) \approx\varphi_2(u_2)$, this means that there exist an element $g = (g_1, \ldots, g_p)\in G^{(p)}$ for some $p\geq n_1, n_2$ and a composition of multiplication maps $\mu_i\colon G^{(p)}\rightarrow G^{(n_i)}$ such that $\mu_i(g) = \varphi_i(u_i)$. There exists a local bisection $g \in \mathbf{c}_i\subset G^{(p)}$ such that $\mu_i(\mathbf{c}_i) = \varphi_i(\mathbf{b}_i)$, because every multiplication map is a submersion. Although $\mathbf{c}_1$ and $\mathbf{c}_2$ might be different, {their points are pairwise Mal'cev equivalent:}
	\begin{equation*}
		AC(\mathbf{c}_1)=AC(\varphi_1(\mathbf{b}_1))=AC(\varphi_2(\mathbf{b}_2))= AC(\mathbf{c}_2),
	\end{equation*}
	wherein we used our assumption on the bisections $\mathbf{b}_1$ and $\mathbf{b}_2$. 
	\[\begin{tikzcd}
		{U_1} && {U_2} \\
		{G^{(n_1)}} & {G^{(p)}} & {G^{(n_2)}} \\
		& {AC(G)}
		\arrow["{\varphi_1}"', from=1-1, to=2-1]
		\arrow["{\varphi_2}", from=1-3, to=2-3]
		\arrow[from=2-1, to=3-2]
		\arrow["{\mu_1}"', from=2-2, to=2-1]
		\arrow["{\mu_2}", from=2-2, to=2-3]
		\arrow[from=2-3, to=3-2]
	\end{tikzcd}\]
			
	For $a = 1, \ldots, p$, let $\{\alpha_{ab}\}_{b = 1}^{q}$ be $q$ elements of $\mathcal{B}$ such that their classes in $\mathcal{B}/I_{t(g_a)}\mathcal{B}$ form a spanning set. The minimal number of $\alpha_{ab}$ needed to form such a spanning set might vary over the different $a$'s, in which case we enlarge the family by adding zero vector fields until the spanning set contains exactly $q$ elements for a fixed $q$. For $i = 1,2$, consider the vector fields 
	\begin{equation*}
		\overrightarrow{\alpha_{ab}}^{(a)}\in \mathfrak{X}(G^{(p)}) \quad\text{ and }\quad d\mu_i(\overrightarrow{\alpha_{ab}}^{(a)})\in \mathfrak{X}(G^{(n_i)}).
	\end{equation*}
	The latter is also of the form $\overrightarrow{\alpha}^{(j)}$ for some $\alpha\in \mathcal{B}$ and $1\leq j \leq n_i$, because the multiplication map $m_k$ between the $k$-th and $(k+1)$-th components maps $\overrightarrow{\alpha}^{(j)}$ to either $\overrightarrow{\alpha}^{(j)}, \overrightarrow{\alpha}^{(j-1)}$ or zero (recall \Cref{eq: multiplication and right-invariant extensions}). Since $U_i$ is a bisubmersion, there exists a lift $X^i_{ab}\in \mathfrak{X}(U_i)$ of $d\mu_i(\overrightarrow{\alpha_{ab}}^{(a)})$. All these vector fields $\{X^i_{ab}\}_{a,b}$ are in $\Gamma(\ker ds_{U_i})$, but they do not need to span the whole of it near $u_i$. Extend them to a spanning set by adding $r$ vector fields  $Y_1^i, \ldots, Y_{r}^i$. The number $r$ might depend on $i$, in which case we add a couple of zero vector fields. Each $Y_c^i$ can be assumed to be $\varphi_i$-related to the zero vector field, because $\varphi_i^{-1}(\overrightarrow{\mathcal{B}}^{(n_i)}) = \Gamma_c(\ker ds_{U_i})$ (see also the proof of Proposition 3.2 in \cite{Zambon})\\
			
	We use the notation $\Phi_{\lambda}^X$ for the time-$\lambda$ flow of $X$ and
	\begin{equation*}
		\Phi_{\lambda_1, \ldots, \lambda_m}^{X_1, \ldots, X_m}:= \Phi_{\lambda_1}^{X_1}\circ \cdots \circ \Phi_{\lambda_m}^{X_m}.
	\end{equation*}
	With this notation, the set $W$ and the maps $p_1,p_2$ can be defined as follows
	\begin{gather*}
		p_i\colon W\subset \left(\mathbf{b}_1\times_M \mathbf{b}_2\right)\times \mathbb{R}^{pq}\times \mathbb{R}^{r}\rightarrow U_i\colon \\
		p_i\left(x^1, x^2, \{\lambda_{ab}\}_{a,b}, \{\lambda_c\}_c\right) = 
		\left(\Phi_{\lambda_1, \ldots, \lambda_{r}}^{Y_1^i, \ldots, Y_{r}^i}\circ\Phi_{\lambda_{11}, \ldots, \lambda_{1q}}^{X_{11}^i, \ldots, X_{1q}^i}\circ \cdots\circ\Phi_{\lambda_{p1}, \ldots, \lambda_{pq}}^{X_{p1}^i, \ldots, X_{pq}^i} \right)\left(x^i\right),
	\end{gather*}
	where $W$ is a neighbourhood of $(u_1, u_2, 0, 0)$ on which the flows (and thus $p_1$ and $p_2$) are	defined. The map $p_i$ is a submersion on a neighbourhood of $u_i$ {because the vectors tangent to the local bisection $\mathbf{b}_i$ provide a complement to} $ \ker ds_{U_i}$, which is generated by $\{X_{ab}^i\}_{a,b}\cup \{Y_c^i\}_c$. Moreover, the element $(u_1, u_2, 0, 0)$ is mapped by $p_i$ to $u_i$. We now show that it makes the following diagram commutative. The arrows and objects coloured in grey do not participate to the commutativity. 
	\[\begin{tikzcd}
		& W \\
		{U_1} && {U_2} \\
		{G^{(n_1)}} & {\textcolor{grey}{G^{(p)}}} & {G^{(n_2)}} \\
		& AC(G)
		\arrow["{p_1}"', from=1-2, to=2-1]
		\arrow["{p_2}", from=1-2, to=2-3]
		\arrow["{\varphi_1}", from=2-1, to=3-1]
		\arrow["{\varphi_2}"', from=2-3, to=3-3]
		\arrow[from=3-1, to=4-2]
		\arrow["{\mu_1}"', color={rgb,255:red,153;green,153;blue,153}, dashed, from=3-2, to=3-1]
		\arrow["{\mu_2}", color={rgb,255:red,153;green,153;blue,153}, dashed, from=3-2, to=3-3]
		\arrow[from=3-3, to=4-2]
	\end{tikzcd}\]
	The diagram commutes if, for a fixed $w = (x^1, x^2, \{\lambda_{ab}\}_{a,b}, \{\lambda_c\}_c)$, the elements $\varphi_1\circ p_1(w)$ and $\varphi_2\circ p_2(w)$ are Mal'cev equivalent.  For $i=1,2$, recall that $X^i_{ab}$ is $\varphi_i$-related to $d\mu_i(\overrightarrow{\alpha_{ab}}^{(a)})$ and $Y_c^i$ to the zero vector field, hence the $n_i$-tuple $\varphi_i\circ p_i(w)$ is Mal'cev equivalent to the $p$-tuple
	\begin{equation}\label{eq: equivalent flows}
		\left(\Phi_{\lambda_{11}, \ldots, \lambda_{1q}}^{\overrightarrow{\alpha_{11}}^{(1)}, \ldots, \overrightarrow{\alpha_{1q}}^{(1)}}\circ \cdots\circ\Phi_{\lambda_{p1}, \ldots, \lambda_{pq}}^{\overrightarrow{\alpha_{p1}}^{(p)}, \ldots, \overrightarrow{\alpha_{pq}}^{(p)}} \right)\left(\mathbf{c}_i(x^i)\right),
	\end{equation}
	where $\mathbf{c}_i(x^i)$ is the unique element on the local bisection $\mathbf{c}_i$ corresponding to $x^i$. {Notice that the flow expression in \Cref{eq: equivalent flows} is independent of $i$; only the initial point $\mathbf{c}_i(x^i)$ depends on $i$.} The bisections $\mathbf{c}_1$ and $\mathbf{c}_2$ were chosen such that $\mathbf{c}_1(x^1)\approx \mathbf{c}_2(x^2)$. {Corollary \ref{cor: flows and malcev} implies that they remain Mal'cev equivalent after applying the flows, hence the diagram commutes}.
\end{proof}

\begin{corollary}\label{cor: morphisms can be locally inverted}
	The relation ``equivalent" of \Cref{def: morphisms using AC(G)} is an equivalence relation on the disjoint union of all bisubmersions.
\end{corollary}
 
Morphisms of bisubmersions can thus be locally inverted in the sense that, if there exists a morphism $U\rightarrow V$ of bisubmersions sending $u\in U$ to $v\in V$, then there also exist an open $V'\subset V$ containing $v$ and a morphism of bisubmersions $V'\rightarrow U$ sending $v$ to $u$.

\subsection{Definition of holonomy groupoid}
	
Given a local Lie groupoid $G$ and a singular subalgebroid $\mathcal{B}$ of $A_G$, we construct a topological groupoid $H(\mathcal{B}, G)$. When the associative completion $AC(G)$ is smooth, we recover the holonomy groupoid $H(\mathcal{B}, AC(G))$ obtained in \cite{Zambon} using the Lie groupoid $AC(G)$, see Proposition \ref{prop: same as in Marco paper when AC(G) is smooth}.

\begin{definition}\label{def: adapted bisubmersion}
	Let $U$ be a bisubmersion and $\mathcal{A} = \{U_i\}_{i\in I}$ a collection of bisubmersions for the same pair $(\mathcal{B}, G)$. We say that $U$ is \textbf{adapted} to $\mathcal{A}$ if for every $u\in U$, there exist an $i\in I$ and a morphism of bisubmersions $f\colon U'\subset U\rightarrow U_i$ defined on a neighbourhood of $u$.
\end{definition}

\begin{definition}\label{def: atlas}
	An \textbf{atlas} $\mathcal{A}$ is a countable collection of bisubmersions for $(\mathcal{B}, G)$ such that 
	\begin{itemize}
		\item $\bigcup_{U\in \mathcal{A}}s_U(U) = M$,
		\item if $U\in \mathcal{A}$, then $U^{-1}$ is adapted to $\mathcal{A}$,
		\item if $U_1, U_2\in \mathcal{A}$, then $U_1\times_M U_2$ is adapted to $\mathcal{A}$. 
	\end{itemize}
\end{definition}

If $\mathcal{U}$ is any countable collection of bisubmersions satisfying the first item of \Cref{def: atlas}, then the collection $\langle \mathcal{U}\rangle$ generated by $\mathcal{U}$, namely all products and inverses of bisubmersions in $\mathcal{U}$, is an atlas. If $\mathcal{U}$ only contains path-holonomy bisubmersions (see \cref{def:phbisub}), then $\langle \mathcal{U}\rangle$ is called a \textbf{path-holonomy atlas}.

\begin{remark}\label{remark: morphisms send path-holonomy bisubemrsion to path-holonomy bisubmersion}
	Let $F\colon G_1\rightarrow G_2$ be a morphism of local Lie groupoids covering $\text{Id}_M$ and $(U,\varphi)$ a path-holonomy bisubmersion for $(\mathcal{B}, G_1)$, then $(U, F\circ \varphi)$ is   a path-holonomy bisubmersion for $(F_*(\mathcal{B}), G_2)$. Indeed, for any $\alpha\in \mathcal{B}_1$, we have
	\begin{equation*}
		F_*(\overrightarrow{\alpha}) = \overrightarrow{F_*(\alpha)}
	\end{equation*}
	thanks to \Cref{lemma: morphisms and right-invariant extensions}. Consequently, the composition
	\begin{gather*}
		F\circ \varphi\colon U\subset \mathbb{R}^n\times M\rightarrow G_2\colon\\
		(\lambda_1, \ldots, \lambda_n, x)\mapsto F\left(\Phi_1^{\sum_i\lambda_i\overrightarrow{\alpha_i}}(x)\right) = \Phi_1^{\sum_i\lambda_i\overrightarrow{F_*(\alpha_i)}}(x)
	\end{gather*}
	is a path-holonomy bisubmersion for $(F_*(\mathcal{B}), G_2)$. Given a path-holonomy atlas $\mathcal{A}$ for $(\mathcal{B}, G_1)$, we thus obtain a path-holonomy atlas for $(F_*(\mathcal{B}), G_2)$ by composing $\varphi$ with $F$ for every path-holonomy bisubmersion $(U, \varphi)\in \mathcal{A}$.
\end{remark}

\begin{definition}\label{def: holonomy groupoid}
	Let $G$ be a local Lie groupoid, $\mathcal{B}$ a singular subalgebroid of $A_G$ and $\mathcal{A}$ a path-holonomy atlas for $(\mathcal{B}, G)$. The \textbf{holonomy groupoid} of the pair $(\mathcal{B}, G)$ is defined as a set by
	\begin{equation}\label{eq: definition of holonomy groupoid}
		H(\mathcal{B} , G):= \Bigg(\bigsqcup_{U\in  \mathcal{A}}U\Bigg) \big/\sim.
	\end{equation}
	where $\sim$ is the equivalence relation that identifies $u_1\in U_1$ with $u_2\in U_2$ if there exists a morphism of bisubmersions $f\colon U_1'\subset U_1\rightarrow U_2$ defined around $u_1$ such that $f(u_1) = u_2$.
\end{definition}
	
\begin{remark}\label{rem:phgroid}
In Definition \ref{def: holonomy groupoid} we defined the  holonomy groupoid choosing a  path-holonomy atlas.
However, we will see in \S \ref{sect: morphism} that different path-holonomy atlases give isomorphic holonomy groupoids.
\end{remark}	

Given an arbitrary (not necessary path-holonomy)  atlas $\mathcal{A}$, we denote by $H(\mathcal{B}, G, \mathcal{A})$  the set defined in \Cref{eq: definition of holonomy groupoid}; it depends on the atlas $\mathcal{A}$. Every result in the sequel up to \Cref{prop: adapted implies morphism between holonomy groupoids} holds for arbitrary atlases, while we specialize to path-holonomy atlases from \Cref{lemma: technical} onward (namely for the functoriality and examples). In particular, \Cref{prop: groupoid structure using AC(G)} shows that $H(\mathcal{B}, G, \mathcal{A})$ is a groupoid for any atlas $\mathcal{A}$, but remark that we reserve the term ``holonomy groupoid" for those constructed out of a path-holonomy atlas.

\begin{theorem}\label{prop: groupoid structure using AC(G)}
	The set $H:=H(\mathcal{B} , G, \mathcal{A})$ is naturally endowed with a topological groupoid structure and a groupoid morphism $\Phi\colon H \rightarrow AC(G)$ covering the identity on $M$. 
\end{theorem}

\begin{proof}
	We show every axiom in the definition of a topological groupoid, see for example \cite[Definition 13.1]{CrainicFernandesMarcut}. The notation $[u\in (U,\varphi)]$ will be used throughout the proof to denote the class of an element $u$ in the domain of  the bisubmersion $\varphi\colon U\rightarrow G^{(n)}$.\\
	
	\underline{The map $\Phi$:}	Define $\Phi$ so that the diagram below commutes, i.e. so that it sends the class $[u\in (U,\varphi)]$ to $(AC\circ \varphi_U)(u)$.
	\[\begin{tikzcd}
		\bigsqcup_{U\in \langle \mathcal{A}\rangle}U & H\\
		\bigsqcup_{n>0}G^{(n)} &   AC(G) \\
		\arrow["/\sim", from=1-1, to=1-2]
		\arrow[ from=1-1, to=2-1]
		\arrow["\Phi", from=1-2, to=2-2]
		\arrow["AC", from=2-1, to=2-2]
	\end{tikzcd}\]
	To show that this map is well-defined, suppose that there exists another element $v \in V$ and a morphism of bisubmersions mapping $u$ to $v$. By definition of a morphism of bisubmersions, we then have $AC(\varphi_U(u)) = AC(\varphi_V(v))$.\\
	 		
	\underline{Source and target:} Define $s_H:= s_{AC(G)} \circ \Phi$ and $t_H:= t_{AC(G)}\circ \Phi$, where $s_{AC(G)}$ is the source of the topological groupoid $AC(G)$ and $t_{AC(G)}$ is its target. Alternatively, the source of $H$ can be defined using the source of a bisubmersion:
	\begin{equation*}
		s_H\left([u \in (U,\varphi)]\right) = s_U(u)
	\end{equation*}
	and similarly for the target. \\
		
	\underline{Multiplication:} Let $u_i \in U_i$ for $i = 1, 2$ such that $s_{U_1}(u_1) = t_{U_2}(u_2)$ and set
	\begin{equation*}
		m_H([u_1], [u_2]):= [(u_1, u_2)\in U_1 \times_M U_2].
	\end{equation*}
	This is well defined, since if $f_i\colon U_i\rightarrow V_i$ is a local morphism of bisubmersions mapping $u_i$ to $v_i$, then 
	\begin{equation*}
		(f_1,f_2)\colon U_1\times_M U_2\rightarrow V_1\times_M V_2
	\end{equation*}
	is a morphism sending $(u_1,u_2)$ to $(v_1,v_2)$. {We briefly verify compatibility of $m_H$ with the source:}
	\begin{align*}
		s_H([(u_1, u_2)\in U_1\times_M U_2]) &= (s_{AC(G)}\circ AC)\left(\varphi_{U_1}(u_1), \varphi_{U_2}(u_2)\right)\\
		 &= (s_{AC(G)}\circ AC)(\varphi_{U_2}(u_2)) = s_H([u_2\in U_2]).
	\end{align*}
	Moreover, associativity of $H$ follows from the fact that taking fibre products is associative, or in symbols $(U_1 \times_M U_2)\times_M U_3 = U_1 \times_M (U_2 \times_M U_3)$. \\
		
	\underline{Units:} Given an arbitrary point $x \in M$, take any $u \in \sqcup_{U\in \mathcal{A}}U$ such that $s_{U}(u) = x$. The unit $x$ inside $H$ is $[(u, u)\in U^{-1} \times_M U]$. We need to show that this association is well-defined. Suppose that both $u$ and $v$ have source $x$ and consider the map
	\begin{equation}\label{eq: multiplication from inside out}
		\begin{aligned}
			M_n\colon \mathcal{M}_n\subset G^{(2n)}&\rightarrow G  \\
			(g_1, \ldots, g_n, h_1, \ldots, h_n)&\mapsto 
			\left(g_1\cdot \left(\cdots \left(\left(g_{n-2}\cdot \left(\left(g_{n-1}\cdot \left(g_n \cdot h_1\right)\right)\cdot h_2\right)\right)\cdot h_3\right)\right)\cdots \right) \cdot h_n.
		\end{aligned}
	\end{equation}
	
	In words, $M_n$ is multiplication from the inside out and its domain $\mathcal{M}_n$ consists of all elements in $G^{(2n)}$ for which $M_n$ is well-defined. Note that $\mathcal{M}_n$ is open and contains all elements of the form $(g_n^{-1}, \ldots, g_1^{-1}, g_1, \ldots, g_n)$. Since $\varphi_{U^{-1}\circ U}(u, u)$ is of this form, the map $M_n\circ \varphi_{U^{-1}\circ U}$ is defined in a neighbourhood of $(u, u)$. {For any local bisection $\mathbf{b}_U$ containing $u$, the image of $(\mathbf{b}_U, \mathbf{b}_U)$ under $\varphi_{U^{-1}\times_MU}$ is Mal'cev equivalent to the unit bisection $s_U(\mathbf{b}_U)$. The same holds with $U$ replaced by $V$, hence} $(u,u)\in U^{-1}\times_M U$ and $(v,v)\in V^{-1}\times_M V$ are equivalent. 
	\[\begin{tikzcd}
		{U^{-1}\times_M U} && {V^{-1}\times_M V} \\
		{G^{(2n)}} & G & {G^{(2m)}} \\
		& {AC(G)}
		\arrow["{(i\circ \varphi_U, \varphi_U)}"', from=1-1, to=2-1]
		\arrow["{(i\circ \varphi_V, \varphi_V)}", from=1-3, to=2-3]
		\arrow["{M_n}", from=2-1, to=2-2]
		\arrow[from=2-1, to=3-2]
		\arrow[from=2-2, to=3-2]
		\arrow["{M_m}"', from=2-3, to=2-2]
		\arrow[from=2-3, to=3-2]
	\end{tikzcd}\]
		
	By construction of the units $M$ inside $H$, it is clear that 
	\begin{gather*}
		s_H(x) = s_H([u\in (U,\varphi)]) = (s\circ \varphi)(u) = x,\\
			t_H(x) = t_H([u\in (U^{-1},i\circ\varphi) ]) = (t \circ i \circ \varphi)(u) = (s\circ \varphi)(u) = x.
	\end{gather*}
	Lastly, we show that $m_H([u], x) = [u]$, for all $[u]$ with $s_H([u]) = x$ and $m_H(x, [v]) = [v]$ for all $[v]$ with $t_H([v]) = x$. For the first one, we have the obvious morphism of bisubmersions
	\begin{equation*}
		U \rightarrow U \times_M U^{-1}\times_M U\colon u\mapsto (u,u,u),
	\end{equation*}
	Notice that the latter represents $m_H([u], x)$. For the second one, note that if $t_H([v\in (V,\varphi_V)]) = x$, then $[v \in (V^{-1},i\circ\varphi_V)]$ satisfies 
	\begin{equation*}
		s_H([v\in V^{-1}]) = (s\circ i \circ \varphi_V)(v) = x
	\end{equation*}
	In particular, the unit $x$ is represented by $(v,v)$ inside $(V^{-1})^{-1}\times_M V^{-1}$. This allows to define the morphism of bisubmersions
	\begin{equation*}
		V\rightarrow V\times_M V^{-1}\times_M V\colon v\mapsto (v,v,v).
	\end{equation*}
	\underline{Inversion:} For an arbitrary element $[u\in (U,\varphi)]$, define its inverse as $i_H([u]) = [u \in (U^{-1},i\circ\varphi)]$. To show that $i_H$ is well-defined, suppose that $u \in U$ and $v \in V$ are related by a morphism of bisubmersions $f$. This means that $AC\circ \varphi_U = AC\circ \varphi_V\circ f$, while the equality we want to show is $AC\circ i\circ \varphi_U = AC\circ i \circ \varphi_V\circ f$ (see the diagram below). The latter does hold, since $AC \circ i = i_{AC} \circ AC$, where $i_{AC}$ is the inversion map of $AC(G)$.
	\[\begin{tikzcd}
		{U^{-1}} && {V^{-1}} \\
		{G^{(n)}} && {G^{(m)}} \\
		& {AC(G)}
		\arrow["f", from=1-1, to=1-3]
		\arrow["{i\circ \varphi_U}"', from=1-1, to=2-1]
		\arrow["{i\circ \varphi_V}", from=1-3, to=2-3]
		\arrow["AC", from=2-1, to=3-2]
		\arrow["AC"', from=2-3, to=3-2]
	\end{tikzcd}\]
		
	By definition of the units, $m_H([u], i_H([u])) = t_H([u])$ and $m_H(i_H([u]), [u]) =s_H([u])$. The equalities $s_H\circ i_H = t_H$ and $t_H \circ i_H = s_H$ follow from the same properties for the inversion $i$ of $G$.\\
		
	Lastly, we show that $\Phi$ is a groupoid morphism covering $\text{Id}_M$. By definition of $s_H$ and $t_H$, it is clear that $\Phi$ intertwines the units, source and target of $H$ and $AC(G)$. Compatibility with multiplication is a direct computation:
	\begin{align*}
		\Phi(m_H([u_1], [u_2])) &= AC(\varphi_{U_1}(u_1), \varphi_{U_2}(u_2))\\
		&= m_{AC}\left((AC\circ\varphi_{U_1})(u_1), (AC\circ\varphi_{U_2})(u_2)\right)\\
		&= 	m_{AC}\left(\Phi([u_1]), \Phi([u_2])\right).
	\end{align*}
\end{proof}

In \Cref{prop: same as in Marco paper when AC(G) is smooth}, we show that the holonomy groupoid defined here is indeed a generalization of the construction by Zambon in \cite{Zambon} whenever $AC(G)$ is a Lie groupoid (in this case, $A_G$ is automatically integrable and $AC(G)$ is one of its integrations). We first need a lemma to relate their atlases of bisubmersions.

\begin{lemma}\label{lemma: bisubmersions from G to AC(G)}
	Let $\varphi\colon U\rightarrow G$ be a simple bisubmersion for $(\mathcal{B} , G)$ and  assume that $AC(G)$ is a Lie groupoid. Then $AC\circ\varphi\colon U\rightarrow AC(G)$ is a bisubmersion for $(\mathcal{B} , AC(G))$.
\end{lemma}
\begin{proof}
	By \cite[Theorem 4.10]{FernandesMichiels}, the map $AC|_G\colon G\rightarrow AC(G)$ is a local diffeomorphism intertwining all the structure maps ($s,t,i$ and $m$) of $G$ and $AC(G)$. The three conditions in \Cref{def: bisubmersion} are now easily verified.
\end{proof}
	
Consequently, any atlas $\mathcal{A} = \langle \{(U_i, \varphi_i)\}_{i\in I}\rangle$ for $(\mathcal{B}, G)$ generated by \textbf{simple} bisubmersions induces an atlas $AC(\mathcal{A}):=\langle \{(U_i, AC\circ \varphi)\}_{i\in I}\rangle$ for $(\mathcal{B}, AC(G))$ by composing each bisubmersion map with $AC$. \Cref{remark: morphisms send path-holonomy bisubemrsion to path-holonomy bisubmersion} for the morphism $AC|_G\colon G\rightarrow AC(G)$ guarantees that  if $\mathcal{A}$ is a path-holonomy atlas, then so is $AC(\mathcal{A})$.

\begin{proposition}\label{prop: same as in Marco paper when AC(G) is smooth}
	If $AC(G)$ is a Lie groupoid and $\mathcal{A}$ is generated by simple bisubmersions, then the groupoids $H_G:= H(\mathcal{B} , G, \mathcal{A})$ of \Cref{def: holonomy groupoid} and $H_{AC(G)}:= H(\mathcal{B} , AC(G), AC(\mathcal{A}))$ from \cite{Zambon} are isomorphic as topological groupoids. 
\end{proposition}
\begin{proof}
	The natural map
	\begin{equation*}
		F\colon H_G\rightarrow H_{AC(G)}\colon [u\in (U,\varphi)]\mapsto[u\in (U,AC\circ\varphi)]
	\end{equation*}
	is well-defined and a morphism of groupoids. It is surjective, since $AC(\mathcal{A})$ contains exactly the bisubmersions $(U, AC\circ \varphi)$ for all $(U, \varphi)\in \mathcal{A}$. For the injectivity, suppose that $F([u_1]) = F([u_2])$ for some $u_i\in U_i$ in the bisubmersion $\varphi_i\colon U_i\rightarrow G^{(n_i)}$. This means that, if we compose both $\varphi_i$ with the completion map $AC$, there exists a local morphism of bisubmersions $f$ making the following diagram commute. 
	\[\begin{tikzcd}
		{U_1'\subset U_1} && U_2 \\
		{G^{(n_1)}} && {G^{(n_2)}} \\
		& {AC(G)}
		\arrow["f", from=1-1, to=1-3]
		\arrow["{\varphi_1}"', from=1-1, to=2-1]
		\arrow["{\varphi_2}", from=1-3, to=2-3]
		\arrow["AC"', from=2-1, to=3-2]
		\arrow["AC", from=2-3, to=3-2]
	\end{tikzcd}\]
	By \Cref{def: morphisms using AC(G)}, the existence of such a local map implies that $u_1$ and $u_2$ are equivalent in $H_G$.
\end{proof}

\subsection{The $\ast$-algebra of the holonomy groupoid}\label{subsec:staralg}

This section is inspired by \cite[\S 4.1-4.3]{AndroulidakisSkandalis} and shows that the $\ast$-algebra construction of a singular subalgebroid given by Androulidakis in \cite[Appendix A]{Zambon} can also be defined when $A$ is not integrable. The discussion here is shorter (and simpler) than in \cite[Appendix A]{Zambon}, since \Cref{prop: equivalence between morphisms and bisections} replaces the three auxiliary lemmas used there.\\

Let $G$ be a local Lie groupoid and $\mathcal{B}$ a singular subalgebroid of $A_G$. Given a bisubmersion $U$, consider the line bundle of complex-valued half densities $|\Lambda|^{1/2}(E_U)$ of the vector bundle $E_U:= \ker ds_U\times \ker dt_U$ and denote its space of compactly supported sections by
\begin{equation*}
	\Omega(U):= \Gamma_c(|\Lambda|^{1/2}(E_U)).
\end{equation*}
We need to understand how these half-densities behave under inverses, multiplications and morphisms of bisubmersions. The following items are proven in \cite[\S 4.2]{AndroulidakisSkandalis} for singular foliations, but the same arguments hold here.

\begin{itemize}
	\item Since $s_{U^{-1}} = t_U$ and $t_{U^{-1}} = s_U$, there is a natural identification between $\Omega(U^{-1})$ and $\Omega(U)$. 
	\item Given two bisubmersions $U_1,U_2$, there is an identification $\Omega(U_1\times_MU_2)\simeq \text{pr}_1^* \Omega(U_1)\otimes_\mathbb{C} \text{pr}_2^* \Omega(U_2)$ using the projection maps $\text{pr}_1\colon U_1\times U_2 \rightarrow U_1$ and $\text{pr}_2\colon U_1\times U_2 \rightarrow U_2$. Given two half-densities $\mu_1\in \Omega(U_1)$ and $\mu_2\in \Omega(U_2)$, there is thus an associated half-density $\mu_1\otimes \mu_2\in\Omega(U_1\times_MU_2)$ defined for every $u_1\in U_1$ and $u_2\in U_2$ via
	\begin{equation*}
		(\mu_1\otimes \mu_2)_{(u_1,u_2)}:= (\mu_1)_{u_1}\otimes (\mu_2)_{u_2}.
	\end{equation*}
	\item If $f$ is a morphism of bisubmersions from $(U, \varphi_U)$ to $(V, \varphi_V)$ that is also a \textbf{submersion}, then $\Omega(U)$ is isomorphic to $ \Gamma_c(|\Lambda|^1(\ker df))\otimes_\mathbb{C} f^*(\Omega(V))$, where $|\Lambda|^1(\ker df)$ stands for the bundle of 1-densities on $\ker df$. This identification gives a \textbf{surjective} map between half densities:
	\begin{equation*}
		f_!\colon \Omega(U)\rightarrow \Omega(V)\colon (f_!(\mu))(x) = \int_{f^{-1}(x)}\mu.
	\end{equation*}
\end{itemize}

The pushforward map $f_!$ is meant to compare half-densities on bisubmersions that are related by a submersive morphism, but not every morphism of bisubmersions is a submersion. However, \Cref{prop: equivalence between morphisms and bisections} says that there exists a third manifold $W$ and two submersions $p\colon W\rightarrow U$, $q\colon W\rightarrow V$. Automatically, both $(W, \varphi_U\circ p)$ and $(W, \varphi_V\circ q)$ become bisubmersions and $p,q$ are (submersive) morphisms of bisubmersions.
\[\begin{tikzcd}
	& W \\
	{U} && {V}
	\arrow["{p}"', from=1-2, to=2-1]
	\arrow["{q}", from=1-2, to=2-3]
\end{tikzcd}\]

\begin{definition}\label{def: star-algebra}
	Let $ \{U_i\}_{i\in I}$ be a path-holonomy atlas for $(\mathcal{B}, G)$ and write $U:= \sqcup_{i\in I}U_i$. A half-density in $\Omega(U)$ is an element of the direct sum $\bigoplus_{i\in I}\Omega(U_i)$. The $\mathbf{\ast}$\textbf{-algebra} of $(\mathcal{B}, G)$ is the quotient
	\begin{equation*}
		\Omega(U)/J = \left(\bigoplus_{i\in I}\Omega(U_i)\right)/J.
	\end{equation*}
	Here $J$ is the ideal generated by those half-densities $\mu\in \Omega(U)$ for which there exists two bisubmersions $W$ and $V$, two submersive morphisms of bisubmersions $p\colon W\rightarrow U$ and $q\colon W\rightarrow V$ and a half-density $\nu \in \Omega(W)$ such that $p_!(\nu) = \mu$ and $q_!(\nu) = 0$. The involution $*$ maps an element $[\mu]\in \Omega(U_i)/J$ to its complex conjugate $[\mu]^*:=[\overline{\mu}]\in \Omega(U_i^{-1})$, but seen as a half-density for the inverse $U^{-1}_i$. The product of two elements $\mu_1\in \Omega(U_1)$ and $\mu_2\in \Omega(U_2)$ reads
	\begin{equation*}
		[\mu_1][\mu_2]:= [\mu_1\otimes \mu_2]\in \Omega(U_1\times_M U_2)/J.
	\end{equation*}
\end{definition}

The proof that both the involution and the product are well-defined is straightforward and similar to \cite[Proposition 4.5]{AndroulidakisSkandalis}; we give them here for completeness. 
\begin{itemize}
	\item \textit{Involution}: Suppose that $\mu\in J$ and let $V,W,p,q,\nu$ be as in the definition of $J$. The inverse bisubmersions $V^{-1},W^{-1}$, the same maps $p,q$ and the complex conjugate $\overline{\nu}\in \Omega(W^{-1})$ satisfy $p_!(\overline{\nu}) = \overline{\mu}$ and $q_!(\overline{\nu}) = \overline{0} = 0$, hence also $\overline{\mu}\in J$.
	\[\begin{tikzcd}[ampersand replacement=\&]
		\& {W^{-1}} \& \\
		{U^{-1}} \&\& {V^{-1}}
		\arrow["p"', from=1-2, to=2-1]
		\arrow["q", from=1-2, to=2-3]
	\end{tikzcd}\]
	\item \textit{Product:} Assume that $\mu_1 \in J$ (the proof is analogous for $\mu_2\in J$) and let $V,W,p,q, \nu$ be as in the definition of $J$. The obvious diagram
	\[\begin{tikzcd}[ampersand replacement=\&]
		\& {W\times_M U_2} \& \\
		{U_1\times_M U_2} \&\& {V\times_M U_2}
		\arrow["{(p, \text{Id}_{U_2})}"'{pos=0.6}, from=1-2, to=2-1]
		\arrow["{(q, \text{Id}_{U_2})}"{pos=0.6}, from=1-2, to=2-3]
	\end{tikzcd}\]
	and the half density $(\nu, \mu_2)\in \Omega(W\times_MU_2)$ show that $\mu_1\otimes \mu_2$ is also in $J$.
\end{itemize} 

\begin{remark}\label{rem:completion}
	The completion of this $*$-algebra to a $C^*$-algebra (either full or reduced) cannot be performed using the exact same procedure as in \cite[\S 4.4 \& 4.5]{AndroulidakisSkandalis}. Indeed, the regular $*$-representations in \cite{AndroulidakisSkandalis} are defined assuming the existence of a holonomy pair, i.e. a pair $(\mathcal{A}, \mathcal{G})$ consisting of an atlas $\mathcal{A}$ and a Lie groupoid $\mathcal{G}$ with smooth source-fibres that is equipped with a morphism of groupoids $\alpha\colon H(\mathcal{F}, \mathcal{A})\rightarrow \mathcal{G}$ such that the composition of $\alpha$ with the projection morphism $U\rightarrow H(\mathcal{F}, \mathcal{A})$ is a leafwise submersion for every $U\in \mathcal{A}$. Debord \cite{Debord2} later showed that such a pair always exists for singular foliations, because the holonomy groupoid itself is leafwise smooth and $\alpha$ can be taken to be the identity map. In our case, the holonomy groupoid $H(\mathcal{B}, G)$ of a singular subalgebroid does not necessarily have smooth source fibres (see \S \ref{sec:Examples}). To find a regular representation of the $*$-algebra given in \Cref{def: star-algebra}, one should give a different description of the $L^2$-space of the source fibre $H(\mathcal{B}, G)_x$ for $x\in M$, which in \cite{AndroulidakisSkandalis} is defined as
	\begin{equation*}
		L^2(H(\mathcal{F})_x):=\overline{\Gamma_c(|\Lambda|^{1/2}(\ker d_xs_H))}.
	\end{equation*}
	A possible solution would be to work at the level of bisubmersions (as in \Cref{def: star-algebra}) and replace the source of the holonomy groupoid $s_H$ by that of a bisubmersion $s_U$ and take the direct sum over all bisubmersions. We will not discuss this here, as the main focus of the current paper is on the geometry underlying singular subalgebroids.
\end{remark}

\section{Morphisms}\label{sect: morphism}

We first show that the holonomy groupoid  $H(\mathcal{B}, G)$ defined in \cref{def: holonomy groupoid} is independent of the choice of  path-holonomy atlas for $\mathcal{B}$, up to isomorphism (see \cref{cor:phiso}). To do so, we show that  certain related atlases and singular subalgebroids induce canonical morphisms between their groupoids. Second, in \S \ref{subsec:imagemor} we establish our main tool to calculate explicitly the holonomy groupoid for a large class of singular subalgebroids, namely those that arise from morphisms of local Lie groupoids (see \Cref{prop: H is powers of K divided by I}). We will apply this tool in \Cref{sec:Examples} to provide non-trivial examples of holonomy groupoids for singular subalgebroids of non-integrable Lie algebroids.\\

The definitions, theorems and proofs are inspired by \cite[\S 3.3.2, \S 4]{Zambon}, but require more care when working with local Lie groupoids. 

\subsection{Adapted atlases}

Let $G$ be a local Lie groupoid and $\mathcal{B}$ a singular subalgebroid of $A_G$. Recall that we defined adapted bisubmersions in \Cref{def: adapted bisubmersion} and atlases in \Cref{def: atlas}.

\begin{definition}
    Suppose that $\mathcal{A}_1, \mathcal{A}_2$ are two atlases for the same pair $(\mathcal{B}, G)$. We say that $\mathcal{A}_1$ is \textbf{adapted} to $\mathcal{A}_2$ if every bisubmersion in $\mathcal{A}_1$ is adapted to $\mathcal{A}_2$.
\end{definition}

Note that if two bisubmersions $U$ and $V$ are adapted to $\mathcal{A}$, then so are the inverse $U^{-1}$ and the product $U\times_M V$. For example, if $f$ is a local morphism from $U$ to some $U'\in  \mathcal{A}$, then the very same map $f$ is also a local morphism from $U^{-1}$ to $U'^{-1}$ (see the following diagram).
\[\begin{tikzcd}[ampersand replacement=\&]
	U \&\& {U'} \\
	{G^{(n_1)}} \& {AC(G)} \& {G^{(n_2)}} \\
	{G^{(n_1)}} \& {AC(G)} \& {G^{(n_2)}}
	\arrow["f", from=1-1, to=1-3]
	\arrow["{\varphi_U}"', from=1-1, to=2-1]
	\arrow["{\varphi_{U_2}}", from=1-3, to=2-3]
	\arrow[from=2-1, to=2-2]
	\arrow["i", from=2-1, to=3-1]
	\arrow["{i_{AC(G)}}", from=2-2, to=3-2]
	\arrow[from=2-3, to=2-2]
	\arrow["i", from=2-3, to=3-3]
	\arrow[from=3-1, to=3-2]
	\arrow[from=3-3, to=3-2]
\end{tikzcd}\]
In particular, if $\mathcal{A}_1$ is generated by a collection of bisubmersions $\{U_i\}_{i\in I}$, then it is enough to check that $U_i$ is adapted to $\mathcal{A}_2$ for every $i\in I$ to conclude that $\mathcal{A}_1$ is adapted to $\mathcal{A}_2$. Recall that the groupoid $H(\mathcal{B}, G, \mathcal{A}_i)$ was defined just before \Cref{prop: groupoid structure using AC(G)}. 

\begin{proposition}\label{prop: adapted implies morphism between holonomy groupoids}
    Let $\mathcal{A}_i$ for $i=1,2$ be two atlases for $(\mathcal{B}, G)$. Write $H_i$ for the groupoid $H(\mathcal{B}, G, \mathcal{A}_i)$ and $\Phi_i\colon H_i\rightarrow AC(G)$ for the canonical map to $AC(G)$. If $\mathcal{A}_1$ is adapted to $\mathcal{A}_2$, then there exists a canonical, injective morphism of groupoids $F\colon H_1\rightarrow H_2$ intertwining $\Phi_1$ and $\Phi_2$. 
\end{proposition}
 
\begin{proof}
    Denote by $\sim_i$ the equivalence relation on $\sqcup_{U\in \mathcal{A}_i}U$ such that $H_i = \sqcup_{U\in \mathcal{A}_i}U/ \sim_i$. Take an arbitrary $u\in U\in \mathcal{A}_1$ and let $f$ be any local morphism of bisubmersions around $u$ with codomain $U_2\in \mathcal{A}_2$. We will show that the map
    \begin{equation*}
    	\tilde{F}\colon\bigsqcup_{U\in \mathcal{A}_1}U\rightarrow H_2\colon u\mapsto [f(u)]
    \end{equation*}
    is independent of the chosen bisubmersion $U_2$ and local morphism $f$. If $f'$ is another such morphism with codomain $V_2$, then \Cref{cor: morphisms can be locally inverted} guarantees the existence of a local morphism $h\colon U_2'\rightarrow U_1$ mapping $f(u)$ to $u$. The composition $g:= f'\circ h$ shows that $f(u)\sim_2 f'(u)$.
    \[\begin{tikzcd}[ampersand replacement=\&]
    	\& {U_2} \\
    	{U_1} \& {U_2'} \\
    	\& {V_2}
    	\arrow["f", from=2-1, to=1-2]
    	\arrow["{f'}"', from=2-1, to=3-2]
    	\arrow[hook, from=2-2, to=1-2]
    	\arrow["h"', from=2-2, to=2-1]
    	\arrow["g", from=2-2, to=3-2]
    \end{tikzcd}\]
	Moreover, $\tilde{F}$ descends to an injective map $F\colon H_1\rightarrow H_2$ after we have shown the following equivalence:
	\begin{equation}\label{eq: proof of adapted implies morphism}
	    \forall u,v\in \bigsqcup_{U\in \mathcal{A}_1}U\colon \left(u\sim_1 v\Leftrightarrow \tilde{F}(u) = \tilde{F}(v)\right).
	\end{equation}
	Suppose that $g$ is a local morphism mapping $u\in U_1$ to $v\in V_1\in \mathcal{A}_1$, then invert $f$ locally in the diagram below to obtain a local morphism $g'$ mapping $f(u)$ to $f'(v)$. 
	\[\begin{tikzcd}
		{U_1} & {U_2} \\
		{V_1} & {V_2}
		\arrow["f", from=1-1, to=1-2]
		\arrow["g"', from=1-1, to=2-1]
		\arrow["{g'}", from=1-2, to=2-2]
		\arrow["{f'}"', from=2-1, to=2-2]
	\end{tikzcd}\]
	Conversely, if $g'$ is a morphism of bisubmersions mapping $f(u)$ to $f'(v)$, then $f'$ can be inverted locally to find a morphism sending $u$ to $v$. This proves \Cref{eq: proof of adapted implies morphism}.\\
	
	The map $F$ is a morphism of groupoids, because if $f\colon U_1\rightarrow U_2$ and $f'\colon V_1\rightarrow V_2$ are morphisms of bisubmersions, then so is 
	\begin{equation*}
		(f,f')\colon U_1\times_MV_1\rightarrow U_2\times_M V_2.
	\end{equation*}
	It remains to show that $\Phi_1 = \Phi_2\circ F$. Choose representatives $u_i\in U_i$ such that $F([u_1]) = [u_2]$, then by construction of $F$, there exists a morphism of bisubmersions $f\colon U_1\rightarrow U_2$ sending $u_1$ to $u_2$. The existence of such a local morphism implies that $\Phi_1([u_1]) = \Phi_2([u_2])$. 
\end{proof}

\begin{lemma}\label{lemma: technical}
	Let $\mathcal{A}$ be a path-holonomy atlas for $(\mathcal{B}, G)$ generated by the set of path-holonomy bisubmersions $\{U_i\}_{i\in I}$. For every $i\in I$, let $V_i\subset U_i$ be {any} neighbourhood of the units $U_i\cap (\{0\}\times M)$. Then $\mathcal{A}$ is adapted to the path-holonomy atlas generated by $\{V_i\}_{i\in I}$.
\end{lemma}

\begin{remark}
	Given a path-holonomy bisubmersion $U$, there is a neighbourhood $V$ of $U\cap(M\times \{0\})$ -- thus $V$ is itself a path-holonomy bisubmersion -- satisfying the property
	\begin{equation}\label{eq: remark smaller path-holonomy bisubmersions}
		s_V(V)=V\cap(M\times \{0\})=t_V(V).
	\end{equation}
	(Furthermore, these sets agree with  $U\cap(M\times \{0\})$).
	Indeed, one can take
	\begin{equation}\label{eq: remark smaller 2}
		V:= s_U^{-1}(U\cap(M\times \{0\}))\cap t_U^{-1}(U\cap(M\times \{0\})).
	\end{equation}
	The properties follow from the fact that $U\cap(M\times \{0\})\subset V$. Therefore, if we take as a new path-holonomy atlas the one generated by the bisubmersions in \Cref{eq: remark smaller 2} as $U$ ranges over all the $U_i$'s, then the original atlas $\mathcal{A}$ is adapted to $\mathcal{A}'$, and all path-holonomy bisubmersions in $\mathcal{A}'$ satisfy \Cref{eq: remark smaller path-holonomy bisubmersions}. 
\end{remark}

\begin{proof}
	Take an arbitrary path-holonomy bisubmersion $(U, \varphi)$ in $\mathcal{A}$ associated to some local generators $\alpha_1, \ldots, \alpha_n\in \mathcal{B}$,  so
	\begin{equation*}
		\varphi(\lambda_1, \ldots, \lambda_n, x) =  \Phi^{\sum_{j = 1}^n\lambda_j\overrightarrow{\alpha}_j}_1(x).
	\end{equation*}
	We want to show that an arbitrary $u = (\lambda_1,\ldots, \lambda_n,x)\in U$ is equivalent to a product of elements inside path-holonomy bisubmersions of the type $V_i$. For brevity, we will denote $\alpha:= \sum_{j = 1}^n\lambda_j\alpha_j\in\mathcal{B}$.
	
	\begin{figure}[h]
		\centering
		\includegraphics[scale=0.7]{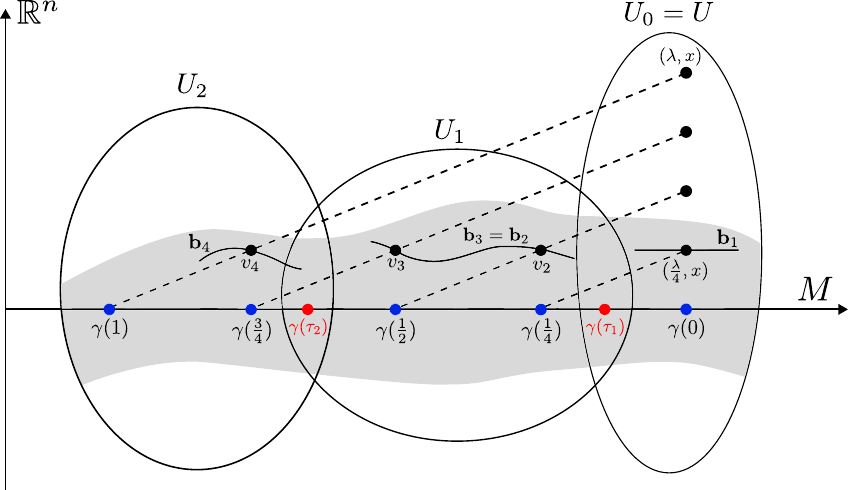}
		\caption{Sketch of the situation for $r = 2$ and $k = 4$. The shaded area represents the union $V_0\cup V_1\cup V_2$.}
	\end{figure}
	
	Consider the path 
	\begin{equation*}
		\gamma\colon [0,1]\rightarrow M\colon \tau \mapsto \left(t\circ \Phi^{\overrightarrow{\alpha}}_{\tau}\right)(x)
	\end{equation*} 
	starting at $x$ and ending at $t(u)$. Every $y\in \Ima \gamma$ is contained in the source of some path-holonomy bisubmersion in $\mathcal{A}$. Since $\gamma$ is compact, we only need finitely many such bisubmersions $U_0, \ldots, U_r$ and we can assume that $\gamma(0)\in U_0$ and $\gamma(1)\in U_r$. We will write $N_i:= s_{U_i}(U_i)$ (they are open subsets of $M$). Since $\gamma$ is connected, we can further rearrange $U_0, \ldots, U_r$ so that $N_i\cap N_{i+1}$ is non-empty for every $i = 0, \ldots, r-1$, and there are $\tau_i\in \gamma^{-1}(N_i\cap N_{i+1})$ such that 
	\begin{equation*}
		\gamma\left(\left[0, \tau_1\right]\right)\subset N_0, \quad \gamma\left(\left[\tau_i, \tau_{i+1}\right]\right)\subset N_i,\quad \gamma\left(\left[\tau_r, 1\right]\right)\subset N_r.
	\end{equation*}
	Let $\{\alpha_{ij}\}_{j = 1}^{n_i}$ be the local generators of $\mathcal{B}$ defining $U_i$. Since $\alpha$ lies in $\mathcal{B}$, we can find functions $f_{ij}\in C^\infty(N_i)$ such that for all $y\in N_i$:
	\begin{equation}\label{eq: change of local generators}
		\alpha_{y} = \sum_{j = 1}^{n_i}f_{ij}(y)(\alpha_{ij})_{y}.
	\end{equation}
	Now take a natural number $k$ large enough such that the following two conditions are satisfied.
	\begin{enumerate}
		\item The compact path
		\begin{equation*}
			[\tau_i,\tau_{i+1}]\rightarrow U_i\colon \tau\mapsto\left(\frac{1}{k}(f_{i1}\circ \gamma)(\tau), \ldots, \frac{1}{k}(f_{in_i}\circ \gamma)(\tau), \gamma(\tau)\right).
		\end{equation*}
		lies completely within $V_i$ for every $i = 1, \ldots, r$.
		\item The compact path
		\begin{equation*}
			[0,1]\rightarrow G^{(2)}\colon \tau \mapsto \left(\Phi_{1/k}^{\overrightarrow{\alpha}}(\gamma(\tau)), \Phi_{\tau}^{\overrightarrow{\alpha}}(x)\right)
		\end{equation*}
		lies completely within the domain of multiplication $G_m$. This is possible because $G_m$ is open and every element of $G$ can be multiplied on the left by its corresponding unit, i.e. $(\gamma(\tau), \Phi_{\tau}^{\overrightarrow{\alpha}}(x))\in G_m$.
	\end{enumerate}
	For every $a = 1, \ldots, k$, define $x_a = \gamma(a/k)$ and let $i(a)$ be any index $i(a) \in \{0,\ldots, r\}$ such that $a/k\in [\tau_{i(a)}, \tau_{i(a)+1}]$, in particular so  that $x_a\in N_{i(a)}$. \\
 
	{\bf Claim:}
	\emph{$u = (\lambda,x)$ is equivalent to the $k$-tuple $(v_k, \ldots, v_1)\in V_{i(k)}\times_M\cdots \times_M V_{i(1)}$, where
	\begin{equation}\label{eq: definition of v_a}
		v_a = \left(\frac{1}{k}f_{i(a)1}(x_a), \ldots, \frac{1}{k}f_{i(a)n_{i(a)}}(x_a), x_a\right).
	\end{equation}}
	Note that $v_a$ lies indeed in $V_{i(a)}$ by item 1. The claim thus shows that $u$ is equivalent to a product $(v_k, \ldots, v_1)$, where each $v_a$ is in some $V_i$ (namely $V_{i(a)}$), which concludes the proof of the lemma.\\

	We now prove the claim. Since $V_{i(a)}$ is open, there is a neighbourhood $N_{i(a)}'\subset N_{i(a)}$ around $x_a$ such that
	\begin{equation*}
		\mathbf{b}_a:= \left\{\left(\frac{1}{k}f_{i(a)1}(y), \ldots, \frac{1}{k}f_{i(a)n_{i(a)}}(y), y\right)\in U_{i(a)}\mid y\in N_{i(a)}'\right\}
	\end{equation*}
	is contained inside $V_{i(a)}$. It is a bisection because \Cref{eq: change of local generators} guarantees that $\mathbf{b}_a$ is a ``constant bisection" (with height $1/k$) w.r.t. $\overrightarrow{\alpha}$. More precisely, if we denote by $\varphi_{i(a)}$ the bisubmersion map of $U_{i(a)}$, then 
	\begin{equation}\label{eq: image of v_a}
		\varphi_{i(a)}\left(\frac{1}{k}f_{i(a)1}(y), \ldots, \frac{1}{k}f_{i(a)n_{i(a)}}(y), y\right) = \Phi_{1/k}^{\sum_{j = 1}^{n_{i(a)}}f_{i(a)j}\overrightarrow{\alpha}_{i(a)j}}(y) = \Phi_{1/k}^{\overrightarrow{\alpha}}(y)
	\end{equation}
	and the right-invariant vector field whose flow is taken (namely $\overrightarrow{\alpha}$) does not depend on the starting point $y$ any more. For the same reason, the constant height section
	\begin{equation*}
		\mathbf{b}:= U\cap \left(\{(\lambda_1, \ldots, \lambda_n)\}\times M\right)
	\end{equation*}
	is a bisection of $U$ containing $u$. Using \Cref{prop: equivalence between morphisms and bisections}, it only rests to show that 
  	\begin{enumerate}
		\item[i)] $\varphi(u)\in G$ is Mal'cev equivalent to $(\varphi_{i(k)}(v_k), \ldots, \varphi_{i(1)}(v_1))$, and
  		\item[ii)] $\varphi(\mathbf{b})$ is Mal'cev equivalent to $\varphi_{i(k)}(\mathbf{b}_{k})\times_M \cdots \times_M \varphi_{i(1)}(\mathbf{b}_{1})$.
	\end{enumerate}
	We will first show i). Recall from \Cref{eq: image of v_a} that $\varphi_{i(a)}(v_a) = \Phi_{1/k}^{\overrightarrow{\alpha}}(x_a)$. Take for example $a = 1$, then item 2. for $\tau = 1/k$ says that $\Phi_{1/k}^{\overrightarrow{\alpha}}(x_1)$ can be multiplied on the right by $\Phi_{1/k}^{\overrightarrow{\alpha}}(x)$ and one finds 
	\begin{equation*}
		\left(\Phi_{1/k}^{\overrightarrow{\alpha}}(x_1), \Phi_{1/k}^{\overrightarrow{\alpha}}(x)\right) \approx R_{\Phi_{1/k}^{\overrightarrow{\alpha}}(x)}\left(\Phi_{1/k}^{\overrightarrow{\alpha}}(x_1)\right) = \Phi_{1/k}^{\overrightarrow{\alpha}}\left(x_1\cdot \Phi_{1/k}^{\overrightarrow{\alpha}}(x)\right) = \Phi_{1/k}^{\overrightarrow{\alpha}}\left(\Phi_{1/k}^{\overrightarrow{\alpha}}(x)\right) = \Phi_{2/k}^{\overrightarrow{\alpha}}(x).
	\end{equation*}
	By induction, if we take $\tau = a/k$, then item 2. guarantees that $\Phi_{1/k}^{\overrightarrow{\alpha}}(x_a)$ can be multiplied by $\Phi_{a/k}^{\overrightarrow{\alpha}}(x)$ and their product is $\Phi_{(a+1)/k}^{\overrightarrow{\alpha}}(x)$. We conclude that the $k$-tuple $(\varphi_{i(k)}(v_k), \ldots, \varphi_{i(1)}(v_1))$ can be multiplied from right to left and that the resulting product is $\Phi_1^{\overrightarrow{\alpha}}(x)$, which is just $\varphi(u)$.\\
	
	The exact same multiplication procedure holds locally around $x$ because the domain $G_m$ is an open subset of $G^{(2)}$. In particular, for an arbitrary $y\in s_U(\mathbf{b})$, if we define $y_a = (t\circ \Phi_{a/k}^{\overrightarrow{\alpha}})(y)$, we obtain that $\varphi(\lambda_1, \ldots, \lambda_n, y)$ is Mal'cev equivalent to the product $(\varphi_{i(k)}(v_k'), \ldots, \varphi_{i(1)}(v_1'))$, where $v_a'$ is as in \Cref{eq: definition of v_a}, but with $x_a$ replaced by $y_a$. Since $N_{i(a)}'$ is open and $x_a\in N_{i(a)'}$, each $v_a'$ will be in $\mathbf{b}_a$ as long as $y$ is close enough to $x$. After making the bisections $\mathbf{b}_a$ smaller if necessary, the image of their product is thus Mal'cev equivalent to that of $\mathbf{b}$. This proves the claim, and thus the lemma.
\end{proof}

\begin{proposition}\label{prop: path is adapted to any atlas}
   If $\mathcal{A}$ is a path-holonomy atlas for $(\mathcal{B}, G)$, then $\mathcal{A}$ is adapted to any atlas of $(\mathcal{B}, G)$. 
\end{proposition}

\begin{proof}
    Write $\mathcal{A}_1:=\mathcal{A}$, let $\mathcal{A}_2$ be any other atlas for $(\mathcal{B}, G)$ and take an arbitrary path-holonomy bisubmersion $U\in \mathcal{A}_1$. Recall that $U\subset \mathbb{R}^n\times M$ is an open neighbourhood of $(0,x)$ for some $x\in M$ with bisubmersion map
    \begin{equation*}
    	\varphi_{U}\colon U\rightarrow G\colon (\lambda_1, \ldots, \lambda_n, y)\mapsto \Phi_1^{\sum_i\lambda_i\overrightarrow{\alpha_i}}(y).
    \end{equation*}
	For an arbitrary $y\in M$, let $U_2\in \mathcal{A}_2$ be any bisubmersion containing an element $u_2$ such that $s_{U_2}(u_2) = y$. For any local bisection $\mathbf{b}\subset U_2$ containing $u_2$, consider the two bisections $(\mathbf{b}, \mathbf{b})\subset U_2^{-1}\times_MU_2$ and $(0, s_{U_2}(\mathbf{b}))\subset U$. They have the same image inside $AC(G)$, namely
    \begin{equation*}
    	\left(AC\circ \varphi_{U_2^{-1}\times_MU_2}\right)(\mathbf{b}, \mathbf{b}) = \left(AC\circ \varphi_{U}\right)(0, s_{U_2}(\mathbf{b})).
    \end{equation*}
    By \Cref{prop: equivalence between morphisms and bisections}, there exists a local morphism from $(0,y)\in U$ to $(u_2,u_2)\in U_2^{-1}\times_MU_2$. This shows that there is an open neighbourhood $V\subset U$ around $U\cap (\{0\}\times M)$ that is adapted to $\mathcal{A}_2$. Doing so for every path-holonomy bisubmersion in $\mathcal{A}_1$, we obtain a new atlas generated by open neighbourhoods $V_i\subset U_i$ for every path-holonomy bisubmersion $U_i$ in $\mathcal{A}_1$. Since the relation ``adapted to" is transitive, \Cref{lemma: technical} finishes the proof. 
\end{proof}

All path-holonomy atlases are thus adapted to each other, by \cref{prop: path is adapted to any atlas}. This and the canonicity of the morphism in \cref{prop: adapted implies morphism between holonomy groupoids} imply the following corollary.
\begin{corollary}\label{cor:phiso}
	Any two path-holonomy atlases for $(\mathcal{B}, G)$ give isomorphic holonomy groupoids.
\end{corollary}
This justifies why we omit the atlas when denoting the holonomy groupoid $H(\mathcal{B}, G)$ in \cref{def: holonomy groupoid}.\\

The holonomy groupoid has a functorial type property stated in the following proposition. Its proof is similar to that of \cite[Proposition 4.1]{Zambon} with minor changes to allow local Lie groupoids.
\begin{proposition}\label{prop: morphisms functioriality}
	Let $F\colon G_1\rightarrow G_2$ be a morphism of local Lie groupoids covering the identity, $\mathcal{B}_1$ a singular subalgebroid of $A_{G_1}$ and $\mathcal{B}_2:=F_*(\mathcal{B}_1)$. Write $H_i:= H(\mathcal{B}_i, G_i)$ for the corresponding holonomy groupoids and $\Phi_i\colon H_i\rightarrow AC(G_i)$ for the canonical map to $AC(G_i)$. There exists a surjective morphism of groupoids $\Xi\colon H_1\rightarrow H_2$ making the following diagram commute.
	\[
	\begin{tikzcd}
		{H_1} & {H_2} \\
		{AC(G_1)} & {AC(G_2)}
		\arrow["\Xi", from=1-1, to=1-2]
		\arrow["{\Phi_1}"', from=1-1, to=2-1]
		\arrow["{\Phi_2}", from=1-2, to=2-2]
		\arrow["{AC(F)}"', from=2-1, to=2-2]
	\end{tikzcd}\]
	If $AC(F)\colon AC(G_1)\rightarrow AC(G_2)$ is injective, then $\Xi$ is an isomorphism.
\end{proposition}

\begin{proof}
	{We already saw in  \Cref{remark: morphisms send path-holonomy bisubemrsion to path-holonomy bisubmersion} that a path-holonomy atlas $\mathcal{A}_1$ for $(\mathcal{B}_1, G_1)$ naturally induces a path-holonomy atlas $\mathcal{A}_2$ for $(\mathcal{B}_2, G_2)$ by composing $\varphi$ with $F$ for every path-holonomy bisubmersion $(U, \varphi)\in \mathcal{A}_1$.} Recall that we are allowed to choose the atlas $\mathcal{A}_2$, because all path-holonomy atlases are adapted to each other. \\
	
	As sets, the disjoint union of all bisubmersions in $\mathcal{A}_1$ is the same as for $\mathcal{A}_2$, but the equivalence relations $\sim_1$ and $\sim_2$ that define $H_1$ and $H_2$ respectively might differ. If there is a local morphism of bisubmersions $f\colon U\rightarrow V$ for some $(U, \varphi_U), (V, \varphi_V)\in  \mathcal{A}_1$, then $f$ is also a morphism of bisubmersions between $(U, F\circ \varphi_U)$ and $(V, F\circ \varphi_V)$, because the following diagram is commutative. 
	\[\begin{tikzcd}
		U && V \\
		{G_1} & {AC(G_1)} & {G_1} \\
		{G_2} && {G_2} \\
		& {AC(G_2)}
		\arrow["f", from=1-1, to=1-3]
		\arrow["{\varphi_U}"', from=1-1, to=2-1]
		\arrow["{\varphi_V}", from=1-3, to=2-3]
		\arrow[from=2-1, to=2-2]
		\arrow["F"', from=2-1, to=3-1]
		\arrow["{AC(F)}"'{pos=0.4}, from=2-2, to=4-2]
		\arrow[from=2-3, to=2-2]
		\arrow["F", from=2-3, to=3-3]
		\arrow[from=3-1, to=4-2]
		\arrow[from=3-3, to=4-2]
	\end{tikzcd}\]
	Herein, we used the functoriality of $AC$ to make the bottom sides of the diagram commute (see \Cref{rem: functoriality of AC} or \cite[\S 4]{FernandesMichiels}). The equivalence classes of $\sim_1$ are thus contained in those of $\sim_2$. The projection $\Xi\colon H_1\rightarrow H_2$ forms a well-defined, surjective morphism of groupoids. The diagram in the statement is commutative by construction. \\
	
	For the last part of the proposition, assume that $AC(F)$ is injective. If $u\sim_2 v$ for two elements $u\in (U, F\circ\varphi_U)$ and $v\in (V, F\circ \varphi_V)$, then there exists a local morphism of bisubmersions {$f\colon U' \rightarrow V$ defined around $u$} such that
	\begin{equation*}
		AC\circ F\circ \varphi_U = AC\circ F\circ \varphi_V\circ f.
	\end{equation*} 
	The functoriality of the $AC$-functor implies $AC\circ F = AC(F)\circ AC$. Since $AC(F)$ is injective, we obtain that $f$ is also a morphism of bisubmersions between $(U, \varphi_U)$ and $(V, \varphi_V)$, namely
	\begin{equation*}
		AC\circ \varphi_U = AC\circ \varphi_V\circ f.
	\end{equation*} 
\end{proof}

\subsection{Singular subalgebroids from morphisms}\label{subsec:imagemor}

Any morphism of local Lie groupoids $\varphi\colon G_1\rightarrow G_2$ covering the identity on the base induces a singular subalgebroid of $A_{G_2}$, namely $\mathcal{B}:= \varphi_*(\Gamma_c(A_{G_1}))$. Such singular subalgebroids are said to \textbf{arise from a morphism of local Lie groupoids} and provide a broad class of examples for which we can calculate the holonomy groupoid $H(\mathcal{B}, G_2)$ explicitly in terms of $G_1$ and $\varphi$. \\

Recall from \Cref{rem: functoriality of AC} that, given such a morphism $\varphi\colon G_1\rightarrow G_2$, there is an induced map $AC(\varphi)$ at the level of associative completions making the diagram in that remark commute.

\begin{definition}
	Let $\varphi\colon G_1\rightarrow G_2$ be a morphism of local Lie groupoid covering the identity. The $\textbf{bisection kernel}$ $\text{Bker}(\varphi)$ of $\varphi$ consists of all elements $[g]\in AC(G_1)$ for which there exist a representative $g\in G_1^{(n)}$ of $[g]$ and a local bisection $\mathbf{b}\subset G_1^{(n)}$ containing $g$ such that $AC(\mathbf{b})\subset \ker(AC(\varphi))$.
\end{definition}

\begin{remark}\label{rem: about definition of bisection kernel}
	A couple of remarks are needed about this definition.
	\begin{enumerate}
		\item The name ``bisection kernel" and its notation are inspired by the terminology ``strong kernel" used in \cite{CamilleRubenLeonidSingularFoliations}, although the meaning is slightly different. Notice that $\text{Bker}(\varphi)$ is contained in the kernel of $AC(\varphi)$ and that it is a normal subgroupoid of $AC(G_1)$.
		
		\item By a local bisection is meant a locally closed submanifold $\mathbf{c}\subset G_1^{(n)}$ such that $s|_\mathbf{c}$ and $t|_\mathbf{c}$ are diffeomorphisms onto open subsets of $M$, where $s$ and $t$ denote the source and target maps of the last, respectively first $G_1$-component.
		
		\item Using the functoriality of $AC$, the composition $AC(\varphi)\circ AC \colon  \sqcup_{n>0}G_1^{(n)}\to AC(G_2)$ can equivalently be written as $AC\circ \varphi$, if we extend the morphism $\varphi$ to a map $\varphi\colon \sqcup_{n>0}G_1^{(n)}\rightarrow \sqcup_{n>0}G_2^{(n)}$ in the obvious way.
		
		\item We consider local bisections on $G^{(n)}$ and not on $AC(G)$, since the latter is not always smooth. However, note that the existence of a local bisection $\mathbf{b}$ satisfying $g\in \mathbf{b}$ and $AC(\mathbf{b})\subset \ker(AC\circ \varphi)$ is independent of the choice of representative  $g  \in G_1^{(n)}$. For example, if $h= m_k(g)\in G^{(n-1)}$ for some $1\leq k\leq n-1$ (i.e. $h$ is a contraction of $g$), then $m_k(\mathbf{b})$ is also a bisection in the kernel of $AC\circ \varphi$. On the other hand, suppose that $m_k(h) = g$ (i.e. $h$ is an expansion of $g$). Using that $m_k$ is a submersion, let $\sigma$ be one of its local sections mapping $g$ to $h$, then $\sigma(\mathbf{b})$ is a local bisection in the kernel of $AC(\varphi)$ containing $h$. Any $h$ in the Mal'cev class of $g$ is obtained via a (finite) series of such contractions and expansions of $g$, hence the belonging of $[g]$ to $\text{Bker}(\varphi)$ is independent of the chosen representative.
	\end{enumerate}
\end{remark}

The bisection kernel of $\varphi$ is a normal subgroupoid of $AC(G_1)$, hence the quotient $AC(G_1)/\text{Bker}(\varphi)$ is a groupoid and the quotient map is a morphism of groupoids.

\begin{proposition}\label{prop: H is powers of K divided by I}
    Let $\varphi\colon G_1\rightarrow G_2$ be a morphism of local Lie groupoids covering the identity and $\mathcal{B}$ a singular subalgebroid of $A_{G_2}$. If $(G_1, \varphi)$ is a bisubmersion for $(\mathcal{B}, G_2)$, then the holonomy groupoid of $(\mathcal{B}, G_2)$ is isomorphic to
    \begin{equation*}
        H(\mathcal{B}, G_2)\cong \frac{AC(G_1)}{\text{Bker}(\varphi)}.
    \end{equation*}
    Moreover, the canonical map $\Phi\colon H(\mathcal{B}, G_2)\rightarrow AC(G_2)$ is the natural map induced by $\varphi$.
\end{proposition}

When $\mathcal{B}= \varphi_*(\Gamma_c(A_{G_1}))$, the assumption above is automatically satisfied, see \cref{prop: morphism is bisubmersion for induced subalgebroid}.

\begin{proof}
    Write $\mathcal{A}$ for a path-holonomy atlas for $(\mathcal{B}, G_2)$ and consider the atlas generated by just one bisubmersion $\mathcal{A}' = \langle\{(G_1, \varphi)\}\rangle$. \Cref{prop: path is adapted to any atlas} says that $\mathcal{A}$ is adapted to $\mathcal{A}'$ and we now show that also $\mathcal{A}'$ is adapted to $\mathcal{A}$.\\
    
    Take an arbitrary unit $x\in M\subset G_1$ and let $(U,\varphi_U)$ be a path-holonomy bisubmersion around $(0,x)\in \mathbb{R}^n\times M$. The bisections $\mathbf{b}_U := U\cap(\{0\}\times M)$ and $\mathbf{b}_{G_1}:= \text{pr}_M(\mathbf{b}_U)$ have equivalent images under their respective bisubmersion maps:
    \begin{equation*}
    	\varphi_U(\mathbf{b}_U)\approx \varphi(\mathbf{b}_{G_1}).
    \end{equation*}
    Using \cref{prop: equivalence between morphisms and bisections}, this means that there exists a neighbourhood $N$ of the units in $G_1$ that is adapted to $\mathcal{A}$. Any neighbourhood of the units generates $G_1$ (\cite[Proposition 2.17]{FernandesMichiels}), meaning that any element in $G_1$ can be written as a well-defined product of elements in $N$. Since a product of adapted bisubmersions is also adapted, we obtain that $G_1$ is adapted to $\mathcal{A}$.\\
    
    Since $\mathcal{A}$ and $\mathcal{A}'$ are adapted to each other, \Cref{prop: adapted implies morphism between holonomy groupoids} implies that their holonomy groupoids are isomorphic. Recall that $\mathcal{A}'$ contains all products and inverses of $G_1$ with itself. The inversion map is a morphism of bisubmersions from $G_1$ to $G_1^{-1}$, hence we do not need to consider inverses:
    \begin{equation*}
    	H(\mathcal{B}, G_2, \mathcal{A}') = \bigsqcup_{n\in \mathbb{N}}G_1^{(n)}/\sim. 
    \end{equation*}
    The equivalence relation $\sim$ identifies two elements $g_i\in G^{(n_i)}$ (for $i = 1,2$) if there exist bisections $\mathbf{b}_i\subset G^{(n_i)}$ containing $g_i$ such that their images under $AC\circ \varphi$ agree. The obvious map
    \begin{equation*}
    	\bigsqcup_{n\in \mathbb{N}}G_1^{(n)}/\sim\ \rightarrow \frac{AC(G_1)}{\text{Bker}(\varphi)}
    \end{equation*}
    is a surjective morphism of groupoids. It is injective by definition of $\text{Bker}(\varphi)$ and $\sim$,  yielding the desired groupoid isomorphism. The statement about $\Phi$ is immediate.
\end{proof}

Proposition 2.13 in \cite{Zambon} remains true for local Lie groupoids. We paraphrase it here for completeness, since it guarantees that the condition in \Cref{prop: H is powers of K divided by I} is automatic when $\mathcal{B}$ arises from $\varphi$.

\begin{proposition}\label{prop: morphism is bisubmersion for induced subalgebroid}
    Let $\varphi\colon G_1\rightarrow G_2$ be a morphism of local Lie groupoids covering the identity and $\mathcal{B}:= \varphi_*(\Gamma_c(A_{G_1}))$ the induced singular subalgebroid. The pair $(G_1, \varphi)$ is a bisubmersion for $(\mathcal{B}, G_2)$. 
\end{proposition}

\begin{proof}
	We check the three conditions of \Cref{def: bisubmersion}. For $i = 1,2$, the source of $G_i$ is denoted $s_i$.
	\begin{enumerate}[I.]
		\item Since $\varphi$ is a morphism of local Lie groupoids, the composition $s_2\circ \varphi = s_1$ is a submersion and similarly for the target.
		\item By assumption, every $\alpha_2\in \mathcal{B}$ is of the form $\varphi_*(\alpha_1)$ for some $\alpha_1\in \Gamma_c(A_{G_1})$. \Cref{lemma: morphisms and right-invariant extensions} exactly states that $\overrightarrow{\alpha_1}$ is $\varphi$-related to $\overrightarrow{\alpha_2}$.
		\item Using \Cref{lemma: pointwise equality is enough for 3rd condition of bisubmersion}, it is enough to show that for all $g\in G_1$, any vector $X\in \ker d_gs_1$ can be extended to a vector field that is $\varphi$-related to $\overrightarrow{\alpha_2}$ for some $\alpha_2\in \mathcal{B}$. Recall that any such vector $X\in \ker d_gs_1$ can be written as $(\overrightarrow{\alpha_1})_g$ for some $\alpha_1\in \Gamma(A_{G_1})$. \Cref{lemma: morphisms and right-invariant extensions} says that $\overrightarrow{\alpha_1}$ is $\varphi$-related to $\overrightarrow{\alpha_2}$ if we set $\alpha_2:= \varphi_*(\alpha_1)$, which concludes the proof.
	\end{enumerate}
\end{proof}

\section{Examples}\label{sec:Examples}

We illustrate how \Cref{prop: H is powers of K divided by I,prop: morphism is bisubmersion for induced subalgebroid} can be used to calculate the holonomy groupoid explicitly in concrete cases. In \S \ref{sec: extreme cases}, we describe the holonomy groupoid in the simple cases that $\mathcal{B} = \Gamma_c(A)$ is as large as possible, or that the base is a point. \Cref{prop: same as in Marco paper when AC(G) is smooth} guarantees that all examples provided in \cite{AndroulidakisSkandalis,Zambon} are also valid examples of the holonomy groupoid constructed here. Consequently, we will focus in \S \ref{sec: quotients of action groupoids} on singular subalgebroids of non-integrable Lie algebroids. Lastly, \S \ref{sec: examples of morphisms between holonomy groupoids} contains examples of the functoriality using \Cref{prop: morphisms functioriality}.

\subsection{General examples}\label{sec: extreme cases}

\begin{example}[The full singular subalgebroid]\label{example: second proof holonomy of largest is AC(G)}
    For any local Lie groupoid $G$, the identity morphism $\varphi\colon G\rightarrow G$ induces the trivial subalgebroid $\mathcal{B} = \Gamma_c(A_G)$ and
    \begin{equation*}
    	H(\mathcal{B}, G) = AC(G).
    \end{equation*}
    Note that this is in contrast with the integrable case, where the holonomy groupoid is smooth if and only if $\mathcal{B}$ is projective, i.e. if and only if $\mathcal{B}\cong \Gamma_c(B)$ for some vector bundle $B$ over $M$ (see \cite[Corollary 4.14]{AndroulidakisZambon}, the case of foliations is discussed in \cite{Debord}). The lack of smoothness of $AC(G)$ can be inherited by the holonomy groupoid.
    {\hfill$\bigtriangleup$}
\end{example}

\begin{remark}[$AC(G_2)$ as a quotient of $AC(G_1)$]\label{remark: ladder-shaped local Lie group}
	Let $\varphi\colon G_1\rightarrow G_2$ be a morphism of local Lie groupoids covering the identity. If $\varphi$ is moreover a submersion, then the induced singular subalgebroid is $\mathcal{B} = \Gamma_c(A_{G_2})$ and hence \Cref{prop: H is powers of K divided by I,example: second proof holonomy of largest is AC(G)} combine to 
	\begin{equation}\label{eq: ladder-shaped}
		AC(G_2)\cong \frac{AC(G_1)}{\text{Bker}(\varphi)}.
	\end{equation}
	This isomorphism can also be obtained directly from $AC(\varphi)\colon AC(G_1)\rightarrow AC(G_2)$, without making use of holonomy groupoids. Indeed, if $\varphi$ is a submersion then it is surjective to a neighbourhood of the units and hence $AC(\varphi)$ is surjective by \cite[Proposition 2.17]{FernandesMichiels}. We now show that the kernel of $AC(\varphi)$ is contained in the bisection kernel of $\varphi$ (the other inclusion is obvious). Suppose that $[g]\in \ker(AC(\varphi))$ for some $g\in G_1^{(n)}$, then there exists a series of contractions and expansions starting at $\varphi(g)$ and ending at $(t(g))\in G_1$ (here we abused notation and wrote $\varphi$ for the map from $G_1^{(n)}$ to $G_2^{(n)}$ naturally induced by $\varphi$). As explained in item 4. of \Cref{rem: about definition of bisection kernel}, there exists a local bisection $\mathbf{c}_2\subset G_2^{(n)}$ containing $\varphi(g)$ that is Mal'cev equivalent to an open neighbourhood of the unit bisection $M$. Since $\varphi$ is a submersion, it admits local sections and thus there exists a local bisection $\mathbf{c}_1\subset G_1^{(n)}$ containing $g$ such that $\varphi(\mathbf{c}_1) = \mathbf{c}_2$. By definition, this means that $[g]$ is in the bisection kernel of $\varphi$.
\end{remark}

\begin{example}\label{ex: ladder-shaped local Lie group}
	As an example of \Cref{eq: ladder-shaped} that will be used later on, consider the open ball $G_1:= B(0, 0.1) $ of radius 1/10 centred at the origin of $(\mathbb{R}^2, +)$ and let $G_2$ be the ladder-shaped local Lie group from \cite[Example 6.1]{FernandesMichiels}. It is known that $AC(G_2)$ is not a smooth manifold, due to  \cite[Theorem 4.10]{FernandesMichiels}. However, \Cref{eq: ladder-shaped} for the inclusion morphism $\varphi\colon G_1\hookrightarrow G_2$ guarantees that $AC(G_2)$ is a quotient of $AC(G_1)\cong \mathbb{R}^2$. One can indeed calculate that 
	\begin{equation*}
		AC(G_2)\cong \frac{\mathbb{R}}{\mathbb{Q}}\times \mathbb{R}.
	\end{equation*}
	{\hfill$\bigtriangleup$}
\end{example}

\begin{remark}[Lie subalgebras]\label{example: local Lie group}
	{The construction of the holonomy group greatly simplifies when the base manifold is a point. Indeed, a single bisubmersion is enough to form an atlas and two elements $u_i\in (U_i, \varphi_i)$ (for $i = 1,2$) are equivalent if and only if $\varphi_1(u_1)\approx \varphi_2(u_2)$.} Let $\mathfrak{g}_2$ be a Lie algebra and $\mathfrak{g}_1\subset \mathfrak{g}_2$ a Lie subalgebra. A global Lie group integrating $\mathfrak{g}_2$ always exists, but let us assume more generally that $G_2$ is a local Lie group with Lie algebra $\mathfrak{g}_2$. Let $G_1$ be any local Lie group integrating $\mathfrak{g}_1$ and $\varphi\colon G_1\rightarrow G_2$ be a morphism integrating the inclusion $\mathfrak{g}_1\hookrightarrow \mathfrak{g}_2$. {Since bisections are just points, \Cref{prop: H is powers of K divided by I} simplifies to}
	\begin{equation}\label{eq: holonomy group arising from morphism}
		H(\mathfrak{g}_1, G_2)\cong \frac{AC(G_1)}{\ker(AC(\varphi))}.
	\end{equation}
	If $G_2$ is a global Lie group, $G_1$ the connected Lie subgroup integrating $\mathfrak{g}_1$ and $\varphi$ the unique embedding integrating the inclusion, we obtain a trivial quotient $H(\mathfrak{g}_1, G_2) = G_1/\{0\} = G_1$. However, the kernel of $AC(\varphi)$ might be non-trivial for general local Lie groups, see \Cref{ex: ladder-shaped local Lie group} where $\mathfrak{g}_1 = \mathfrak{g}_2 = \mathbb{R}^2$.
\end{remark}

\begin{example}\label{example: subalgebras of ladder-shaped}
	To illustrate \Cref{eq: holonomy group arising from morphism}, consider the ladder-shaped local Lie group $G_2$ from \cite{FernandesMichiels} that appeared in \Cref{ex: ladder-shaped local Lie group}. Its Lie algebra $\mathfrak{g}_2$ is abelian; let $\mathfrak{g}_1\subset \mathfrak{g}_2$ be any one-dimensional Lie subalgebra. For $G_1$, take the local Lie subgroup of the open ball of radius 1/10 inside $(\mathbb{R}^2, +)$ integrating $\mathfrak{g}_1$ with the obvious inclusion $\varphi\colon G_1\rightarrow G_2$. The holonomy group is always $H(\mathfrak{g}_1, G_2) = \mathbb{R}$, since $\ker(AC(\varphi))$ is trivial, except when $\mathfrak{g}_1 = \mathbb{R}\oplus \{0\}$, in which case
	\begin{equation*}
		H(\mathfrak{g}_1, G_2) = \frac{\mathbb{R}}{\mathbb{Q}}.
	\end{equation*}
	Although all one-dimensional Lie subalgebras are related by an automorphism of $\mathfrak{g}_2$, not every one of them integrates to an automorphism of $G_2$; otherwise their holonomy groups would be isomorphic. This is a consequence of the fact that Lie algebroid morphisms are equivalent to germs of local Lie groupoid morphisms \cite{CabreraMarcutSalazar}, but not every such morphism globalizes to the whole local Lie groupoid.
	{\hfill$\bigtriangleup$}
	\end{example}	
	
\subsection{Quotients of action groupoids}\label{sec: quotients of action groupoids}
 
\begin{remark}
Let $\mathfrak{g}$ be a Lie algebra, $A$ a Lie algebroid over $M$, and $\psi\colon \mathfrak{g}\to \Gamma(A)$ a Lie algebra map.  Then $\mathcal{B}=\langle \psi(v)\mid v\in \mathfrak{g} \rangle$ is a singular subalgebroid of $A$.
 Composing $\psi$ with the anchor map, we obtain an infinitesimal action of $\mathfrak{g}$ on $M$, hence an action Lie algebroid $\mathfrak{g}\ltimes M$. The Lie algebra map $\psi$ is encoded by the Lie algebroid morphism $$\mathfrak{g}\ltimes M\to A, (v,x)\mapsto \psi(v)_x$$ covering $\text{Id}_M$, whose image at the level of compactly supported sections is precisely $\mathcal{B}$. Hence the holonomy groupoid of $\mathcal{B}$ can be described as a quotient
 of an action groupoid, by \Cref{prop: H is powers of K divided by I,prop: morphism is bisubmersion for induced subalgebroid}. The following three examples arise in this way, taking $\mathfrak{g}$ to be the one-dimensional Lie algebra.
\end{remark} 

\begin{example}[Maps into the local Lie groupoid of \Cref{example: local Lie groupoid over S^2}, I]\label{ex: zero vector field}
	Consider the Lie algebroid $TM\oplus_{\omega} \mathbb{R}_M$ over $M = S^2\times S^2$, with bracket twisted by the 2-form $\omega = \text{pr}_1^*\omega_{S^2}+\lambda\text{pr}_2^*\omega_{S^2}\in \Omega^2(M)$, for $\lambda$ a real number {and  $\omega_{S^2}$ the canonical volume form on the $2$-sphere (with total volume $1$)}. For any $f\in C^{\infty}(M)$,  consider the singular subalgebroid $\mathcal{B}$ generated by the section $(0,f)$. In this example we exhibit its holonomy groupoid $H(\mathcal{B}, G_2)$, where $G_2$ is the local Lie groupoid of \Cref{example: local Lie groupoid over S^2}. Recall that $G_2$ is a subset of $M\times M\times \mathbb{R}$ and that its associators are uniform, but not discrete:
	\begin{equation*}
		\text{Assoc}(G_2) = \Delta M\times \Lambda,
	\end{equation*}
	where $\Delta M$ is the diagonal of $M\times M$ and $\Lambda\subset\mathbb{R}$ is a countable subgroup. Let $G_1 = M\times \mathbb{R}$ be the trivial bundle of Lie groups. Let  $f\in C^\infty(M)$ and consider 
	the morphism
	\begin{equation*}
		\varphi\colon G_1\rightarrow G_2\colon (x,a)\mapsto (x,x, f(x) a).
	\end{equation*}
	We will use \cref{prop: H is powers of K divided by I} to calculate the holonomy groupoid of $\mathcal{B}:= \varphi_*(\Gamma_c(A_{G_1}))$. 	Note that the induced morphism of Lie algebroids $d\varphi\colon A_{G_1}\rightarrow A_{G_2}$ was studied in \cite[Example 2.7]{AndroulidakisAntonini}.\\
	
	In this case, $AC(G_1)$ is naturally identified with $G_1$ and $AC\circ \varphi = \varphi$. A local bisection around $(x, a)\in G_1$ is the same as a function $g\in C^\infty(N)$ defined on a neighbourhood $N\subset M$ of $x$ such that $g(x) = a$. Hence $(x,a)\in \text{Bker}(\varphi)$ if and only if there exists such a local function $g$ satisfying
	\begin{equation}\label{eq: example of local Lie groupoid morphism with non-integrbale algebroid}
		\forall y\in N\colon f(y)g(y)\in \Lambda.
	\end{equation}
	We distinguish three cases.
	\begin{itemize}
		\item If $f$ vanishes on a neighbourhood of $x$, then any $g$ defined on that neighbourhood satisfies \Cref{eq: example of local Lie groupoid morphism with non-integrbale algebroid}. Hence $\text{Bker}(\varphi)_x$  is the whole fibre $\mathbb{R}$ at $x$.
		\item If $f(x)\neq 0$, then there exists a neighbourhood $N\subset M$ of $x$ where $f$ does not vanish. A necessary condition for \Cref{eq: example of local Lie groupoid morphism with non-integrbale algebroid} is that $f(x)a\in \Lambda$ (plug in $y = x$). This is actually a sufficient condition, as one sees taking the local function
		\begin{equation*}
			g\colon N\rightarrow \mathbb{R}\colon y\mapsto \frac{f(x)a}{f(y)}
		\end{equation*}
		Hence $(x,a)\in \text{Bker}(\varphi)$ if and only if $f(x)a\in \Lambda$ .
		\item If $f(x) = 0$, but every neighbourhood of $x$ contains a point where $f$ does not vanish, then \Cref{eq: example of local Lie groupoid morphism with non-integrbale algebroid} cannot be satisfied (except when $a = 0$) for any $g\in C^\infty(M)$ defined around $x$. Indeed, after   making it smaller if necessary, $N$ can be assumed to be connected and $g$ nowhere vanishing {(since we assume $a\neq 0$)}. The product $fg\in C^\infty(N)$ vanishes at $x$, but not at some $y\in N$. By the intermediate value theorem, it attains every value between zero and $f(y)g(y)$. {This contradicts \Cref{eq: example of local Lie groupoid morphism with non-integrbale algebroid}, since $\Lambda$ is countable.} \Cref{eq: example of local Lie groupoid morphism with non-integrbale algebroid} is thus false for every local function $g$ and hence $(x,a)\notin \text{Bker}(\varphi)$. {We conclude that $\text{Bker}(\varphi)_x= \{0\}$.}
	\end{itemize}
	The holonomy groupoid of $(\mathcal{B}, G_2)$ is the quotient of $G_1 = M\times \mathbb{R}$ for which the source fibre at $x\in M$ is
	\begin{equation*}
		H(\mathcal{B}, G_2)_x = 
		\begin{cases}
			\{0\} & \text{ if $f$ vanishes identically around $x$,}\\
			\frac{\mathbb{R}}{\Lambda/f(x)} & \text{ if $f(x)\neq 0$,}\\
			\mathbb{R} & \text{ else}.
		\end{cases}
	\end{equation*}
	Although $\mathcal{B}$ is relatively simple, its holonomy groupoid is already quite pathological.  
	{\hfill$\bigtriangleup$}
\end{example}
	
\begin{example}[Maps into the local Lie groupoid of \Cref{example: local Lie groupoid over S^2 times R}]\label{example: map onto local Lie groupoid over S2 times R}
	Recall from \cite{FernandesMichiels} that a Lie algebroid is integrable if and only if it integrates to a local Lie groupoid with uniformly discrete associators. In \Cref{ex: zero vector field}, the associators were uniform, but non-discrete in every fibre. The current example is very similar, but the associators are not uniform.\\

	Let $M = S^2\times I$ be the product of the sphere with  an open interval containing zero, and denote by $a$ the coordinate on $I$. Consider the Poisson bivector $\pi:= (1+a^2)\omega_{S^2}^{-1}\in \mathfrak{X}^2(M)$ and the associated Lie algebroid $T^*M$. For any $f\in C^{\infty}(S^2)$, let $\mathcal{B}$ be the singular subalgebroid generated by the section $f\,da$. In this example, we exhibit its holonomy groupoid $H(\mathcal{B}, G_2)$, where $G_2$ is the local Lie groupoid of \Cref{example: local Lie groupoid over S^2 times R}. Recall that $G_2$ is a subset of $S^2\times S^2\times I\times \mathbb{R}$ with associators at $(x,a)\in M$ given by 
	\begin{equation*}
		\text{Assoc}_{(x,a)}(G_2) = \{(x,x,a)\}\times \Lambda_a
	\end{equation*}
	for $\Lambda_a\subset\mathbb{R}$ a countable subgroup containing $\phi(a)\mathbb{Z}$, where  $\phi(a) = \frac{-2a}{(1+a^2)^2}$. Let $G_1 = M\times \mathbb{R}$ be the trivial bundle of Lie groups. Given a function $f\in C^\infty(S^2)$, we will calculate the holonomy groupoid of the singular subalgebroid $\mathcal{B}$ of $T^*M=A_{G_2}$ obtained from the morphism
 	\begin{equation*}
		\varphi\colon G_1 \rightarrow G_2\colon (x,a,b)\mapsto (x,x,a,f(x) b).
	\end{equation*}
	Similar to \Cref{ex: zero vector field}, $\text{Bker}(\varphi)$ consists of those elements $(x,a,b)\in G_1$ for which there exists a function $g\in C^\infty(N)$ defined around $(x,a)$ such that $g(x,a)=b$ and
	\begin{equation}\label{eq: bisection example over S^2 times R}
		\forall (x',a')\in N\colon f(x')g(x',a')\in \Lambda_{a'}.
	\end{equation}
	Whether such a $g$ exists depends on the vanishing of $a$ and $f(x)$.
	
	\begin{itemize}
		\item If $f$ vanishes on a neighbourhood of $(x,a)$, then \Cref{eq: bisection example over S^2 times R} is trivially satisfied for any $g$. {Hence $\text{Bker}(\varphi)_{(x,a)}$ is the whole fibre $\mathbb{R}$.}
		\item If $f(x) = 0$, but $f$ is not identically zero around $x$, then $(x,a,b)\in\text{Bker}(\varphi)$ if and only if $b = 0$. {Indeed, suppose that $g\in C^\infty(N)$ would satisfy \Cref{eq: bisection example over S^2 times R} and assume without loss of generality that $N$ is connected and $g$ nowhere vanishing. Take $(x',a')\in M$ so that $(x',a')\in N$ and $f(x')\neq 0$, then $f(x)g(x,a') = 0$, while $f(x')g(x',a')\neq 0$. If we fix $a'$ and consider the function $y\mapsto f(y)g(y,a')$, then the intermediate value theorem forces this function to attain every value between $0$ and $f(x')g(x',a')$. This contradicts the fact that $\Lambda_{a'}$ is countable, hence such a $g$ cannot exist.}
		\item Assume that $f(x)\neq 0$ and $a=0$. Since $\Lambda_0 = \{0\}$ is trivial, \Cref{eq: bisection example over S^2 times R} immediately forces $b = 0$ (plug in $x' = x$ and $a' = a$). In that case, $g$ can be chosen to vanish identically. The bisection kernel at $(x,a)$ is thus trivial.
		\item If $f(x)\neq 0 \neq a$, then a necessary condition for the existence of a $g$ satisfying \Cref{eq: bisection example over S^2 times R} is $f(x) b\in \Lambda_a$. The latter is also sufficient: let $N$ be a neighbourhood where $f(x')\neq 0 \neq a'$ for all $(x',a')\in N$, then
		\begin{equation*}
			g\colon N\rightarrow \mathbb{R}\colon (x',a')\mapsto \frac{\phi(a')}{\phi(a)}\frac{f(x)}{f(x')}b
		\end{equation*}
		satisfies $g(x, a) = b$ and \Cref{eq: bisection example over S^2 times R} thanks to \Cref{eq: relation between associators of S2 times I}. The bisection kernel at $(x,a)$ is
		\begin{equation*}
			\text{Bker}(\varphi)_{(x,a)} = \frac{\Lambda_a}{f(x)}.
		\end{equation*}
	\end{itemize}
	We conclude that $H(\mathcal{B}, G_2)$ is the quotient of $M\times \mathbb{R}$ with source fibre at $(x,a)\in M$ given by
	\begin{equation*}
		H(\mathcal{B},G)_{(x,a)} =
		\begin{cases}
			\{0\} & \text{ if $f$ vanishes identically around $x$,}\\
			\cfrac{\mathbb{R}}{\Lambda_a/f(x)} & \text{ if $f(x)\neq 0\neq a$,}\\
			\mathbb{R} & \text{ else.}
		\end{cases}
	\end{equation*}
	{\hfill$\bigtriangleup$}
\end{example}

The following example is a variation of \Cref{ex: zero vector field}, in which we consider the singular subalgebroid of $TM\oplus_{\omega} \mathbb{R}_M$ generated by $(X,f)$ for a specific choice of vector field $X$ on $M=S^2\times S^2$ and invariant function $f$.

\begin{example}[Maps into the local Lie groupoid of \Cref{example: local Lie groupoid over S^2}, II]\label{example: transformation Lie algebroid}
	
	Consider the Lie algebroid $TM\oplus_{\omega} \mathbb{R}_M$ over $M = S^2\times S^2$, with bracket twisted by the 2-form $\omega = \text{pr}_1^*\omega_{S^2}+\lambda\text{pr}_2^*\omega_{S^2}\in \Omega^2(M)$, for $\lambda$ a real number. Let $X\in \mathfrak{X}(M)$ be the vector field
	\begin{equation*}
		X = \frac{\partial}{\partial \theta_1}+ \frac{\partial}{\partial \theta_2},
	\end{equation*}
	where $\theta_i$ is the polar angle in cylindrical coordinates on the $i$-th copy of $S^2$. In other words, $X$ is the vector field of simultaneous rotations along the vertical axes with period $2\pi$. Let $f\in C^\infty(M)$ be a function that is constant along the integral curves of $X$:
	\begin{equation}\label{eq: function constant along vector field flow}
		\forall a\in \mathbb{R}\colon f\circ\Phi_a^X = f, 
	\end{equation}
	namely $f(x_1, x_2)$ only depends on the relative angle $\theta_2-\theta_1$ and the heights $h(x_i)$ of $x_i\in S^2$ with respect to the axis of rotation. In this example, we exhibit the holonomy groupoid $H(\mathcal{B}, G_2)$, where $\mathcal{B}$ is (globally) generated by the section $(X,f)$ and $G_2$ is the local Lie groupoid of \Cref{example: local Lie groupoid over S^2}. The singular subalgebroid $\mathcal{B}$ can be obtained from a morphism of local Lie groupoids as follows. Let $M\ltimes \mathbb{R}$ be the action groupoid of the infinitesimal action 
	\begin{equation*}
		\mathbb{R}\rightarrow \mathfrak{X}(M)\colon a\mapsto aX.
	\end{equation*}
	Its Lie algebroid $A_{G_1}$ is the vector bundle $\mathbb{R}_M := M\times \mathbb{R}$ with vanishing bracket {on constant sections} and anchor $(x, a)\mapsto aX_x$. The morphism
	\begin{equation}\label{eq: morphism in last example}
		\mathbb{R}_M\rightarrow TM\oplus_\omega \mathbb{R}_M\colon (x,a)\mapsto (aX_x, f(x) a)
	\end{equation}
	at the level of sections induces $\mathcal{B} = \langle (X,f)\rangle$ and it must integrate to a morphism of local Lie groupoids
	\begin{equation*}
		\varphi\colon G_1\subset M \ltimes \mathbb{R}\rightarrow G_2 \subset M\times M\times \mathbb{R}.
	\end{equation*}
	\emph{Step 1: description of the morphism $\varphi$.}	The first two components must be given by $(x,a)\mapsto (\Phi_a^X(x), x)$, which only lies inside $G_2$ when $a$ is not an odd multiple of $\pi$. Hence, the domain of $\varphi$ is the shrinking of $M\ltimes \mathbb{R}$ to the open subset
	\begin{equation*}
		G_1:= \{(x,a)\in M\ltimes \mathbb{R}\mid a\notin (2\mathbb{Z}+1)\pi\}.
	\end{equation*}
	To define the third component of $\varphi$, we introduce an {auxiliary} function $f_a\in C^\infty(S^2)$ depending on $a$. Write $k:= \lfloor a/2\pi\rfloor \in \mathbb{Z}$, then $f_a$ assigns to a point on the sphere $p\in S^2$ the value $$f_a(p)=k\pi+  \sphericalangle(p, a),$$ where $\sphericalangle(p, a)\in [0,\pi]$ is the angle at $p$ between the rotation vector field $V:= \frac{\partial}{\partial \theta}$ and the geodesic joining $p$ with $\Phi_a^{V}(p)$, see \cref{fig: geodesic triangle on sphere}. Note that this angle only depends on the height of $p$ with respect to the rotation axis, which implies that it is invariant under the flow of $V$:
	\begin{equation}\label{eq: f_a invariant under rotation flow}
		\forall a,b\in \mathbb{R}\colon f_b\circ \Phi_{a}^V = f_b.
	\end{equation}
	
	\begin{figure}[h]
		\centering
		\includegraphics[scale = 0.45]{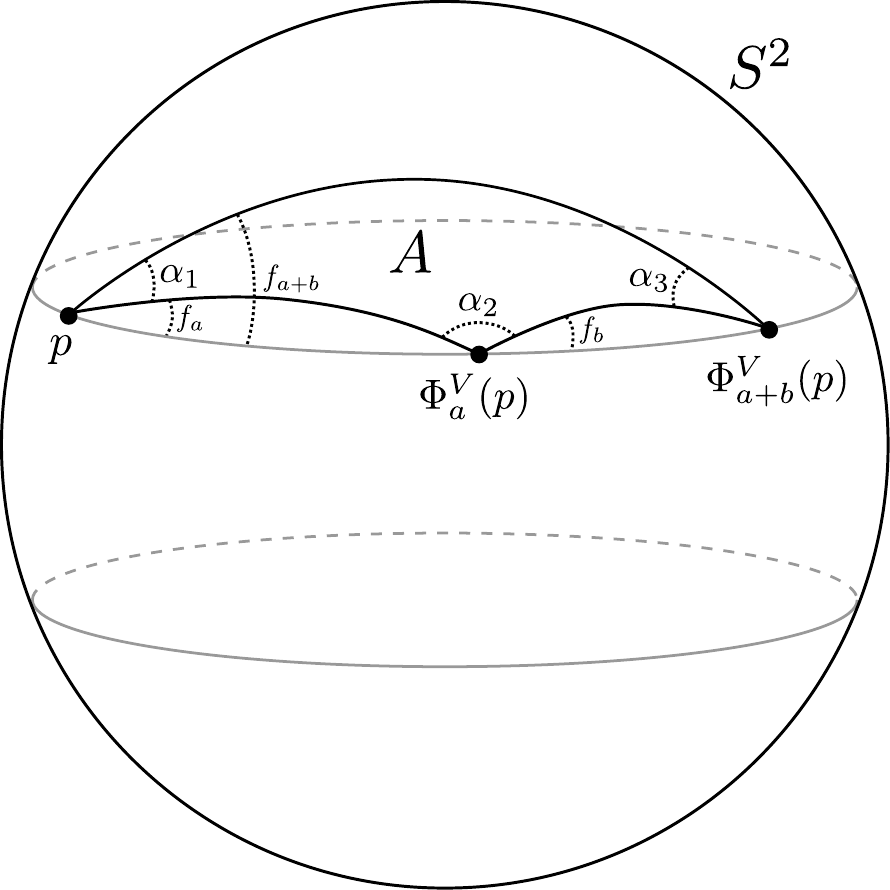}
		\caption{Geodesic triangle on $S^2$ formed by $p, \Phi_a^{V}(p)$ and $\Phi_{a+b}^{V}(p)$}
		\label{fig: geodesic triangle on sphere}
	\end{figure}
	
	We claim that the third component of the morphism $\varphi$ is the following
	\begin{gather*}
		\varphi_3\colon G_1\subset S^2\times S^2\times \mathbb{R}\rightarrow \mathbb{R}\colon (x_1,x_2,a)\mapsto\\
		f(x)a+2\left(f_a(x_1)-a\frac{d}{d\varepsilon}\Bigr|_0f_\varepsilon(x_1)\right)+2\lambda\left(f_a(x_2)-a\frac{d}{d\varepsilon}\Bigr|_0f_\varepsilon(x_2)\right).
	\end{gather*}
 	Indeed, the derivative of this map at $a = 0$ sends $(x_1, x_2,a)$ to $f(x)a$, which is exactly the last component of \Cref{eq: morphism in last example}. To show that $\varphi$ is a morphism of groupoids, recall {from \Cref{example: local Lie groupoid over S^2}} that the multiplication of $G_2$ is 
	\begin{equation}\label{eq: formula area 1}
		(z,y,a)(y,x,b) = (z,x,a+b+A(z_1,y_1,x_1)+\lambda A(z_2,y_2,x_2)),
	\end{equation}
	where $x = (x_1, x_2)$, $y = (y_1, y_2)$, $z = (z_1, z_2)$ and  $A(z_i,y_i,x_i)$ is the signed area of the geodesic triangle formed by the three points on the sphere $z_i,y_i$ and $x_i$. For the triple $x_i, y_i =\Phi_a^{V}(x_i)$ and $z_i = \Phi_{a+b}^{V}(x_i)$, this area can be calculated using the Gauss-Bonnet theorem (see \cref{fig: geodesic triangle on sphere}) as
	\begin{equation}\label{eq: formula area 2}
	 	\begin{split}
	 		A(\Phi_{a+b}^{V},(x_i),\Phi_{a}^{V}(x_i),x_i) &= \alpha_1+\alpha_2+\alpha_3-\pi\\
	 		& = 2(f_{a+b}(x_i)-f_a(x_i)-f_b(x_i)).
	 	\end{split}
	\end{equation}
	Using the invariance assumption \eqref{eq: function constant along vector field flow} on the function $f$ and \Cref{eq: f_a invariant under rotation flow,eq: formula area 1,eq: formula area 2}, it is a straightforward calculation to check that $\varphi$ respects the multiplications of $G_1$ and $G_2$.\\
	
	\emph{{Step 2: description of $\text{Bker}(\varphi)$.}}	\Cref{prop: H is powers of K divided by I} says that $H(\mathcal{B}, G_2)$ is isomorphic to the quotient of $AC(G_1) = M\ltimes \mathbb{R}$ by $\text{Bker}(\varphi)$. We claim that $\text{Bker}(\varphi) = M\times \{0\}$ is trivial and hence that $$H(\mathcal{B}, G_2)\cong M\ltimes \mathbb{R}.$$ The units are always in $\text{Bker}(\varphi)$, therefore we take an arbitrary $(x, a)\in M\ltimes \mathbb{R}$ with $a\neq 0$ and show that it is not in the bisection kernel. \\
	
	Note that $\text{Bker}(\varphi)\subset G_1(\subset AC(G_1))$, since $G_1$ contains all isotropic elements of $M\ltimes \mathbb{R}$. Hence $(x, a)$ is in $\text{Bker}(\varphi)$ if and only if there exists a local function $g\in C^\infty(N)$ on a neighbourhood of $x$ with $g(x) = a$ and $\varphi(y, g(y))\in \text{Assoc}(G_2)$ for all $y\in N$.
	\[\begin{tikzcd}
		{G_1} & {G_2} \\
		{M\ltimes \mathbb{R}} & {AC(G_2)}
		\arrow["\varphi", from=1-1, to=1-2]
		\arrow[from=1-1, to=2-1]
		\arrow[from=1-2, to=2-2]
		\arrow["{AC(\varphi)}"', from=2-1, to=2-2]
	\end{tikzcd}\]
	Recall that the associators of $G_2\subset M\times M\times \mathbb{R}$ are given by $\Delta M\times \Lambda$, {where $\Lambda$ is the countable subgroup of $\mathbb{R}$ defined in 	\Cref{example: local Lie groupoid over S^2}}. If they contain $\varphi(y, g(y))$, then the flow $\Phi_{g(y)}^X$ must fix $y$ and hence $g(y) = 2\pi k$ for some non-zero $k\in \mathbb{R}$ (if $k = 0$, then $a = 0$). To calculate the value $\varphi_3(y,2\pi k)$, we need a more explicit formula for $f_\varepsilon$. {Using the slerp formula \cite{Shoemake} to find the geodesic joining two points on $S^2$ and basic Riemannian geometry}, one obtains
	\begin{equation*}
		\cos(f_\varepsilon(y_i)) = \frac{\sqrt{1-h(y_i)^2}\sin(\varepsilon)}{\sqrt{1-((1-h(y_i)^2)\cos(\varepsilon)+h(y_i)^2)^2}},
	\end{equation*}
	where $y = (y_1, y_2)\in M$ and $h(y_i)$ is the height of $y_i$ w.r.t. the axis of rotation. The derivative of $f_\varepsilon(y_i)$ at $\varepsilon = 0$ is $\frac{h(y_i)}{2}$ and hence
	\begin{equation*}
		\varphi_3(y_1,y_2, 2\pi k) = 2\pi k(2+\lambda - h(y_1)-\lambda h(y_2)).
	\end{equation*}
	Note that this expression depends continuously on the heights $h(y_1)$ and $h(y_2)$, {leading to a contradiction}. Indeed, we assumed that $\varphi_3(y, g(y))$ is in $\Lambda$ for every $y\in N$, but $\Lambda$ {is countable}. We conclude that the local bisection $g$ cannot exist for any $(x,a)\in G_1$ that is not a unit and hence the bisection kernel is trivial.
	{\hfill$\bigtriangleup$}
\end{example}

\subsection{Morphisms between holonomy groupoids}\label{sec: examples of morphisms between holonomy groupoids}

We now turn to examples of \Cref{prop: morphisms functioriality}. The morphism $\Xi\colon H(\mathcal{B}, G_1)\rightarrow H(F_*(\mathcal{B}), G_2)$ induced by $F\colon G_1\rightarrow G_2$ is always surjective, but it does not need to be injective in general, even if $F$ is injective.
\begin{example}[An inclusion of local Lie groups]\label{example: inclusion of local Lie groups}
	Take as in \Cref{ex: ladder-shaped local Lie group} the open ball $G_1 := B(0, 0.1)\subset (\mathbb{R}^2, +)$ with the induced local Lie group structure. Let $G_2$ be the ladder-shaped, non-associative local Lie group from \cite[Example 6.1]{FernandesMichiels} and $F\colon G_1\rightarrow G_2$ the inclusion morphism. The induced map between the holonomy groupoids of $\mathfrak{g}_1:=A_{G_1}$ and $\mathfrak{g}_2:= F_*(\mathfrak{g}_1)$ is the quotient map
	\begin{equation}\label{eq: inclusion of local Lie groups}
		\Xi\colon {H(\mathcal{B}_1, G_1)} \cong \mathbb{R}^2\rightarrow {H(\mathcal{B}_2, G_2)} \cong \mathbb{R}\times \frac{\mathbb{R}}{\mathbb{Q}}.
	\end{equation}
	\Cref{eq: inclusion of local Lie groups} illustrates that our construction of the holonomy groupoid is not intrinsic. {It does depend on the ambient Lie algebroid $A$ in which the singular subalgebroid $\mathcal{B}$ lies}.
	{\hfill$\bigtriangleup$}
\end{example}

\begin{example}\label{example: isotropic foliation}
	The holonomy groupoid of a singular foliation $\mathcal{F}$ constructed in \cite{AndroulidakisSkandalis} uses the Lie algebroid $TM$ with integration the pair groupoid $G_2:= M\times M$. In this example, we {relate} $H(\mathcal{F}, G_2)$ with the holonomy groupoid of the same foliation seen as a singular subalgebroid of $TM\oplus_\omega\mathbb{R}_M$, {where $\omega$ is a closed 2-form on $M$}. Note that $\mathcal{F}$ is  involutive inside $TM\oplus_\omega\mathbb{R}_M$ only when it is isotropic for $\omega$. {This is trivially satisfied for the foliation given below, since its rank is one.}\\
	
	Concretely, let  $G_1\subset  M\times M\times \mathbb{R}$ be the local Lie groupoid of \Cref{example: local Lie groupoid over S^2} over $M = S^2\times S^2$ with Lie algebroid $TM\oplus_\omega\mathbb{R}_M$. {The projection   $\text{pr}\colon G_1\rightarrow G_2$ is a morphism of local Lie groupoids, integrating the projection $TM\oplus_\omega\mathbb{R}_M\to TM$.} Consider the foliation $\mathcal{F}\subset \mathfrak{X}(M)$ globally generated by $X = \frac{\partial}{\partial \theta_1}+\frac{\partial}{\partial \theta_2}$, the vector field of simultaneous rotation that already appeared in \Cref{example: transformation Lie algebroid}. The usual holonomy groupoid of this foliation is the transformation groupoid $H(\mathcal{F}, M\times M) = M\ltimes S^1$. If $\mathcal{F}$ is viewed as a singular subalgebroid of $TM\oplus_\omega\mathbb{R}_M$, then $H(\mathcal{F}, G_1) = M\ltimes \mathbb{R}$, as already calculated in \Cref{example: transformation Lie algebroid} using $f = 0$. The morphism of \Cref{prop: morphisms functioriality} {between holonomy groupoids} is also the obvious projection:
	\begin{equation*}
		\Xi\colon M\ltimes \mathbb{R}\rightarrow M\ltimes S^1.
	\end{equation*}
	To see this, remark that the single path-holonomy bisubmersion 
	\begin{gather*}
		\varphi\colon U:= \{(x, a)\in M\times \mathbb{R}\mid a\notin (2\mathbb{Z}+1)\pi\}\rightarrow G_1\colon\\
		(x,a)\mapsto \Phi_a^{\overrightarrow{X}}(x,x,0) 
	\end{gather*}
	generates an atlas for $(\mathcal{F}, G_1)$, while $(U, \text{pr}\circ \varphi)$ forms an atlas for $(\mathcal{F}, G_2)$. For all integers $k\in \mathbb{Z}$, there exists an obvious morphism of bisubmersions
	\begin{equation*}
		f_k\colon (U, \text{pr}\circ \varphi)\rightarrow (U, \text{pr}\circ \varphi)\colon (x,a)\mapsto (x,a+2\pi k),
	\end{equation*}
	hence $(x, a)$ and $(x,a+2\pi k)$ are identified  inside the usual holonomy groupoid $H(\mathcal{F}, G_2)$. However, they are not identified inside $H(\mathcal{F}, G_1)$, as was proven at the end of \Cref{example: transformation Lie algebroid}, see the beginning of Step 2 there.
	{\hfill$\bigtriangleup$}
\end{example}

\bibliographystyle{amsalpha}

\providecommand{\bysame}{\leavevmode\hbox to3em{\hrulefill}\thinspace}
\providecommand{\MR}{\relax\ifhmode\unskip\space\fi MR }
\providecommand{\MRhref}[2]{%
  \href{http://www.ams.org/mathscinet-getitem?mr=#1}{#2}
}
\providecommand{\href}[2]{#2}

\end{document}